\documentclass[11pt,epsfig]{article}
\usepackage{amsmath}
\usepackage{amsthm}
 \usepackage{amssymb}
\usepackage{epsfig}
\usepackage{leftidx}
\usepackage{graphicx}
\usepackage{amsfonts}
\usepackage{graphicx, color, epstopdf}
\usepackage{cite}
\usepackage{enumerate}
\usepackage{enumitem}
\usepackage{multicol}
\usepackage{booktabs}
\usepackage [latin1]{inputenc}
\usepackage{multirow}
\usepackage{algpseudocode,algorithm,algorithmicx}

\newtheorem{theorem}{Theorem}[section]
\newtheorem{lemma}{Lemma}[section]

\newtheorem{example}{Example}

\newcommand{\vparl}[2]{\ensuremath{\tfrac{\delta #1}{\delta #2}}}

\newcommand{\diff}{\triangledown_{\!\tau}}
\newcommand{\defeq}{:=}
\newcommand{\zd}{\,\mathrm{d}}
\newcommand{\abs}[1]{\left|#1\right|}

\newcommand{\bra}[1]{\left(#1\right)}
\newcommand{\brab}[1]{\big(#1\big)}
\newcommand{\braB}[1]{\Big(#1\Big)}
\newcommand{\brat}[1]{(#1)}
\newcommand{\kbra}[1]{\left[#1\right]}
\newcommand{\kbrab}[1]{\big[#1\big]}

\newcommand{\mynorm}[1]{\left\|#1\right\|}
\newcommand{\mynormb}[1]{\big\|#1\big\|}
\newcommand{\mynormB}[1]{\Big\|#1\Big\|}

\title{An energy stable and maximum bound preserving scheme
with variable time steps for time fractional Allen-Cahn equation\thanks{Updated on \today}}
\author{Hong-lin Liao\thanks{ORCID 0000-0003-0777-6832;
Department of Mathematics,
Nanjing University of Aeronautics and Astronautics,
Nanjing 211106, P. R. China. E-mails: liaohl@nuaa.edu.cn and liaohl@csrc.ac.cn.
This author's work is supported by NSF of China under grant number 12071216.}
\quad Tao Tang\thanks{Division of Science and Technology, BNU-HKBU United International College,
    Zhuhai, Guangdong Province, and Department of Mathematics and International Center for Mathematics, Southern
    University of Science and Technology, Shenzhen, Guangdong Province,China.
    Email: tangt@sustech.edu.cn. This author's work is partially supported by NSF of China under grant numbers 11731006 and K20911001.}
\quad Tao Zhou\thanks{LSEC, Institute of Computational Mathematics and Scientific/Engineering Computing,
Academy of Mathematics and Systems Science, Chinese Academy of Sciences, Beijing, 100190, China. Email: tzhou@lsec.cc.ac.cn. This author's work is partially supported by NSF of China (under grant numbers 11822111 and 11688101), science challenge project (No. TZ2018001), and Youth Innovation Promotion Association of CAS.}
}
\date{}

\begin{document}

\maketitle

\begin{abstract}
In this work, we propose a Crank-Nicolson-type scheme with variable steps for the time fractional Allen-Cahn equation. The proposed scheme is shown to be unconditionally stable (in a variational energy sense), and is maximum bound preserving. Interestingly, the discrete energy stability result obtained in this paper can recover the classical energy dissipation law when the fractional order $\alpha \rightarrow 1.$ That is, our scheme can asymptotically preserve the energy dissipation law in the $\alpha \rightarrow 1$ limit. This seems to be the first work on variable time-stepping scheme that can preserve both the energy stability and the maximum bound principle.

Our Crank-Nicolson scheme is build upon a reformulated problem associated with
the Riemann-Liouville derivative. As a by product, we build up a reversible transformation between the L1-type formula of the Riemann-Liouville derivative and a new L1-type formula of the Caputo derivative, with the help of a class of discrete orthogonal convolution kernels. This is the first time such a \textit{discrete} transformation is established between two discrete fractional derivatives.
We finally present several numerical examples with an adaptive time-stepping strategy to show the effectiveness of the proposed scheme.
\end{abstract}


\section{Introduction}
\setcounter{equation}{0}

In this work, we are concerned with numerical methods for the following time fractional Allen-Cahn (TFAC) equation,
\begin{align}\label{Cont: TFAC}
\partial_t^\alpha u = \varepsilon^2\Delta u - f(u)\quad\text{for $x\in\Omega$ and $t>0$}.
\end{align}
Here $\Omega=(0,L)^2,$ and $\partial_t^\alpha:={}_0^C{\mathcal{D}}_t^{\alpha}$ is the Caputo derivative of order $\alpha$,
\begin{align}\label{Cont: Caputo def}
\partial_t^\alpha v={}_0^C{\mathcal{D}}_t^{\alpha}v\defeq \mathcal{I}_t^{1-\alpha}v'\quad\text{for $0<\alpha<1$},
\end{align}
where $\mathcal{I}_t^{\mu}$ is the Riemann-Liouville fractional integration operator of order $\mu > 0$
\begin{align}\label{Cont: RL integral def}
\bra{\mathcal{I}_t^{\mu}v}(t)\defeq\int_{0}^{t}\omega_{\mu}(t-s)v(s)\zd{s},\quad
\textmd{with} \quad  \omega_{\mu}(t)\defeq t^{\mu-1}/\Gamma(\mu).
\end{align}
The nonlinear bulk force $f=F'(u)$ is given by
\begin{align}
F(u)=\frac14(1-u^2)^2.
\end{align}
For simplicity, we consider periodic solution $u$ along the boundary.

The above time fractional Allen-Cahn equation has been studied both theoretically and numerically in recent years \cite{IncYusufAliyuBaleanu:2018Time,DuYangZhou:2019,LiWangYang:2017,TangYuZhou:2018On,JiLiaoZhang:2020,ZhaoChenWang:2019On,LiuChengWangZhao:2018Time,QTY}. When $\alpha \rightarrow 1,$ the TFAC equation recovers the classical Allen-Cahn equation \cite{AllenCahn_1979}:
\begin{align}\label{Cont: AC}
\partial_t u = \varepsilon^2\Delta u - f(u).
\end{align}
Note that equation (\ref{Cont: AC}) is an $L^2$ gradient flow of a free energy, i.e.,
\begin{align}\label{Cont: free energy}
\partial_t u :=-\frac{\delta E}{\delta u} \quad\text{where the energy}\quad E[u](t) = \int_\Omega \left( \frac{\varepsilon^2}{2} |\nabla u|^2 + F(u)\right) \textmd{dx}.
\end{align}
In this sense, one may view the TFAC equation as a \textit{fractional} gradient flow:
\begin{align}
\partial_t^\alpha u :=-\frac{\delta E}{\delta u}.
\end{align}
It is well known that for the classical AC equation (\ref{Cont: AC}), there holds the energy dissipation law
\begin{align}\label{Cont: EDL}
\frac{\zd E}{\zd t}+\mynormb{\vparl{E}{u}}^2=0, \quad \textmd{or} \quad E[u](t)\leq{E}[u](s), \quad \forall t>s,
\end{align}
and the maximum bound principle
\begin{align}\label{Cont: MBP}
|u(\mathbf{x},t)|\leq{1} \quad \textmd{if} \quad  |u(\mathbf{x},0)|\leq 1.
\end{align}

Thus it is natural to ask whether the TFAC equation (\ref{Cont: TFAC}) also preserves these two properties.
In \cite{TangYuZhou:2018On}, it was shown that the TFAC equation also admits the maximum bound principle (\ref{Cont: MBP}).  However, one can only obtain the following energy stability \cite{TangYuZhou:2018On}:
\begin{align}\label{Cont: ES}
E[u](t)\leq{E}[u](0).
\end{align}
Notice that this is different from the energy dissipation law (\ref{Cont: EDL}).

While it may be interesting to further check whether (\ref{Cont: EDL}) holds for TFAC equation,  however, as a new fractional gradient flow, we shall investigate in this work  a new (yet natural) energy law for the TFAC equation.

\subsection{A variational energy dissipation law}

The first aim of this work is to define a new variational energy dissipation law. To this end, we first rewrite the TFAC into
an equivalent form that involves the Riemann-Liouville
derivative.

Recall the Riemann-Liouville derivative ${}^R\!\partial_t^{\alpha}:={}_{\;\;\,0}^{RL}{\mathcal{D}}_t^{\alpha}$ defined by
\begin{align}\label{Cont: Riemann-Liouville def}
{}^R\!\partial_t^{\alpha} v\defeq \partial_t\mathcal{I}_t^{1-\alpha}v,  \quad\text{for $0<\alpha<1$}.
\end{align}
Due to the semigroup property of the fractional integral we have
\begin{align}\label{Cont: RLderivative property1}
{}^R\!\partial_t^{1-\alpha}\bra{\partial_t^\alpha v}=\partial_t\mathcal{I}_t^{1}v'=v'\quad\text{for $0<\alpha<1$},
\end{align}
Thus, one can reformulate the TFAC equation \eqref{Cont: TFAC} into the following form
\begin{align}\label{Cont: AllenCahn RL}
\partial_tu =-{}^R\!\partial_t^{1-\alpha}\bra{\vparl{E}{u}}.
\end{align}
Moreover, for the Riemann-Liouville derivative of order $1-\alpha$ there holds \cite{AlsaediAhmadKirane:2015}
\begin{align}\label{RL property}
v(t)\brab{{}^R\!\partial_t^{1-\alpha}v}(t)\geq\frac12\brab{{}^R\!\partial_t^{1-\alpha}v^2}(t)+\frac12\omega_{\alpha}(t)v^2(t), \quad \forall v\in C[0,T].
\end{align}
Now, we take the inner product of \eqref{Cont: AllenCahn RL} by $\vparl{E}{u}$ to obtain
\begin{align}\label{Weak formula}
\frac{\zd E}{\zd t}=\brab{\partial_tu,-\varepsilon^2\Delta u+f(u)}
&= -\braB{\vparl{E}{u},{}^R\!\partial_t^{1-\alpha}\vparl{E}{u}},
\end{align}
where $\bra{\cdot,\cdot}$ denotes the $L^2$ inner product.
The above discussions motivate us to define the following variational energy functional $\mathcal{E}_{\alpha}$:
\begin{align}\label{Cont: pseudo-local energy}
\mathcal{E}_{\alpha}[u]:=E[u]+\frac12\mathcal{I}_t^{\alpha}\mynormb{\vparl{E}{u}}^2.
\end{align}
Then, by (\ref{RL property}) and (\ref{Weak formula}) it is easy to show that for $\mathcal{E}_{\alpha}$ it holds
\begin{align}\label{new law}
\frac{\zd \mathcal{E}_{\alpha}}{\zd t}+\frac12\omega_{\alpha}(t)\mynormb{\vparl{E}{u}}^2\leq0, \quad \forall t>0.
\end{align}
That is, the functional $\mathcal{E}_{\alpha}$ seems to be a naturally defined variational energy that admits the dissipation law.  More importantly, when the fractional order $\alpha\rightarrow1$, the above energy law recovers the classical energy dissipation law of AC equation, i.e.,
\begin{align*}
\frac{\zd E}{\zd t}+\mynorm{\vparl{E}{u}}^2\leq0, \quad \forall t>0.
\end{align*}
In this sense, definition (\ref{Cont: pseudo-local energy}) is asymptotically energy dissipation preserving in the $\alpha\rightarrow1$ limit.

\newpage
\subsection{Summary of main contributions}

Our main contribution is two folds:
\begin{itemize}
\item  We design a Crank-Nicolson-type scheme with variable steps for the TFAC equation that can preserve
the new variational energy law \eqref{new law}. The proposed scheme is also shown to
preserve the maximum bound \eqref{Cont: MBP}. Moreover, the discrete variational energy stability can also recover the classical discrete energy dissipation law
when the fractional order $\alpha \rightarrow 1.$ In other words, at the discrete level our scheme can asymptotically preserve the energy dissipation law. This seems to be the first work on variable time-stepping scheme that can preserve both the energy stability and the maximum bound principle.

\item The proposed Crank-Nicolson scheme is build upon the reformulated problem \eqref{Cont: AllenCahn RL} that involves the Riemann-Liouville derivative \eqref{Cont: Riemann-Liouville def}. As a by product, we build up a reversible transformation between the L1-type formula of the Riemann-Liouville derivative \eqref{Cont: Riemann-Liouville def} and a new L1-type formula of the Caputo derivative \eqref{Cont: Caputo def}, with the help of a class of discrete orthogonal convolution kernels. This is the first time such a \textit{discrete} transformation is established between the two discrete fractional derivatives.

\end{itemize}
Finally, we present several numerical examples with an adaptive time-stepping strategy to show the effectiveness of the proposed scheme.

The rest of the paper is organized as follows. In Section 2, we present our numerical scheme and show the discrete variational energy dissipation law. Section 3 is devoted to the unique solvability of our scheme and discrete maximum bound principle. This is followed by some numerical examples in Section 4. We finally give some concluding remarks in Section 5.

\section{Numerical schemes}
\setcounter{equation}{0}
This section will be devoted to the design of our structure preserving Crank-Nicolson type scheme. All our discussions will be emphasized on nonuniform time grids, and this is motivated by the fact that nonuniform grids are powerful in capturing the multi-scale behaviors (including the singular behavior near the initial time) for time-fractional Allen-Cahn equation.

To begin, we consider the following nonuniform time grids:
$$0=t_{0}<t_{1}<\cdots<t_{k-1}<t_{k}<\cdots<t_{N}=T$$
with the step sizes $\tau_{k}:=t_{k}-t_{k-1}$ for $1\leq{k}\leq{N}$.
Let the maximum time-step size $\tau:=\max_{1\leq{k}\leq{N}}\tau_{k}$ and the adjoint time-step ratios $r_n:=\tau_n/\tau_{n-1}$ for $n\ge2$.
Always, we assume the summation $\sum_{k=i}^{j}\cdot=0$ and the product $\prod_{k=i}^{j}\cdot=1$ for index $i>j$.

\subsection{Discrete Riemann-Liouville derivative}
Our scheme will be designed upon the equivalent form (\ref{Cont: AllenCahn RL}).
Consider a mesh function $v^{k}=v(t_k),$  we set (for $k\geq{1}$)
$$\triangledown_{\tau}v^{k}:=v^{k}-v^{k-1}, \quad \partial_{\tau}v^{k-\frac12}:=\triangledown_{\tau}v^{k}/\tau_k, \quad
v^{k-\frac{1}{2}}:=(v^{k}+v^{k-1})/2.$$
Let $(\Pi_{0,k}v)(t)$ be the constant interplant of a function $v(t)$ at $t_{k-1}$ and $t_{k}$, then a piecewise constant approximation is defined as
\begin{align}\label{eq: constants interpolation}
\Pi_{0}v:=\Pi_{0,k}v\quad
\text{so that}\quad
(\Pi_{0,k}v)(t)=v^{k-\frac12}\quad
\text{for $t_{k-1}<{t}\leq t_{k}$ and $k\geq1$}.
\end{align}
For any fixed $n\ge1$, we consider the following discrete Riemann-Liouville derivative for \eqref{Cont: Riemann-Liouville def},
\begin{align}\label{eq: L1r formula}
({}^R\!\partial_{\tau}^{1-\alpha}v)^{n-\frac{1}{2}}
:=&\,\frac{1}{\tau_{n}}\int_{t_{n-1}}^{t_{n}}
\frac{\partial}{\partial t}\int_{0}^{t}\omega_{\alpha}(t-s)(\Pi_{0}v)(s)\zd{s}\zd{t}
\triangleq\frac{1}{\tau_{n}}\sum_{k=1}^{n}a_{n-k}^{(n)}v^{k-\frac{1}2}
\end{align}
for $n\geq1$. The associated discrete convolution kernels $a_{n-k}^{(n)}$ are defined as follows
\begin{align}\label{eq: L1r discrete kernels}
a_{0}^{(n)}:=q_{0}^{(n)}>0\;\; \text{for $n\ge 1$}\quad\text{and}\quad
a_{n-k}^{(n)}:=q_{n-k}^{(n)}-q_{n-k-1}^{(n-1)}<0\quad \text{for $n\ge k+1\ge2$},
\end{align}
where we have used the following auxiliary sequence
\begin{align}\label{eq: DCO kernels for L1r}
q_{n-k}^{(n)}:=\int_{t_{k-1}}^{t_k}\omega_{\alpha}(t_n-s)\zd{s}=\sum_{j=k}^na_{j-k}^{(j)}>0\quad\text{for $1\leq k\leq n$.}
\end{align}

The numerical approximation formula \eqref{eq: L1r formula} has been investigated in \cite{Mustapha:2011,MustaphaAlMutawa:2012}
for linear subdiffusion problems, and the approximation order is shown to be $1+\alpha.$ Notice that this formula was originally called the L1 formula of the Riemann-Liouville derivative \eqref{Cont: Riemann-Liouville def}. However, to avoid possible confuses, we call it here L1$_R$ formula to distinguish it from another well-known L1 formula \cite{LiaoMcLeanZhang:2019,LiaoYanZhang:2019} of the Caputo derivative \eqref{Cont: Caputo def}.

As shown in \cite[Section 2]{MustaphaMcLean:2009},
the kernel of the the Riemann-Liouville derivative ${}^R\!\partial_t^{1-\alpha}$ is positive semi-definite, i.e.,
\begin{align}\label{ieq: positive Riemann-Liouville}
\int_0^Tv(t)\brab{{}^R\!\partial_t^{1-\alpha} v}(t)\zd t
=\int_{0}^Tv(t)\frac{\partial}{\partial t}\int_0^t\omega_{\alpha}(t-s)v(s)\zd s\zd{t}\ge0\quad\text{for $v\in L^2[0,T]$}.
\end{align}
The above L1$_R$ formula \eqref{eq: L1r formula} is designed in a structure preserving way, more precisely, we have
\begin{align*}
\sum_{j=1}^n\tau_jv^{j-\frac12}({}^R\!\partial_{\tau}^{1-\alpha}v)^{j-\frac{1}{2}}
=&\,\sum_{j=1}^n\tau_j\brab{\Pi_{0,j}v}({}^R\!\partial_{\tau}^{1-\alpha}v)^{j-\frac{1}{2}}\nonumber\\
=&\,\int_{t_{0}}^{t_{n}}(\Pi_{0}v)(t)\frac{\partial}{\partial t}\int_{0}^{t}\omega_{\alpha}(t-s)(\Pi_{0}v)(s)\zd{s}\zd{t}
\ge0\quad\text{for $n\ge1$.}
\end{align*}
The arbitrariness of function $v$ implies that the discrete L1$_R$ kernels $a_{n-k}^{(n)}$ in \eqref{eq: L1r discrete kernels} are positive semi-definite. As shown in \cite{Mustapha:2011,MustaphaAlMutawa:2012}, this property implies the $L^2$ norm stability
of numerical scheme when the L1$_R$ formula is applied to linear diffusion $\partial_tu ={}^R\!\partial_t^{1-\alpha}\Delta{u}+f$.
Nevertheless, we remark that the numerical analysis in this work is new and quite different from those in \cite{Mustapha:2011,MustaphaAlMutawa:2012}
as we have to deal with the nonlinear term.

In the next, we show that the discrete kernels $a_{n-k}^{(n)}$ are positive definite without using the
continuous property \eqref{ieq: positive Riemann-Liouville}. This result will be used to establish
the discrete variational energy dissipation law in the forthcoming sections.
\begin{lemma}\label{lem: L1r positive definite}
For any real sequence $\{w_k\}_{k=1}^n,$ the discrete convolution kernels $a_{n-k}^{(n)}$  and $q_{n-k}^{(n)}$ defined
in \eqref{eq: L1r discrete kernels}-\eqref{eq: DCO kernels for L1r} satisfy
\begin{align*}
2w_k\sum_{j=1}^ka_{k-j}^{(k)}w_j\ge
w_k^2\sum_{j=1}^ka_{k-j}^{(k)}+\sum_{j=1}^kq_{k-j}^{(k)}w_j^2-\sum_{j=1}^{k-1}q_{k-j-1}^{(k-1)}w_j^2\quad \text{for $k\ge1$}
\end{align*}
 so that the discrete kernels $a_{n-k}^{(n)}$ are positive definite in the sense that
\begin{align*}
2\sum_{k=1}^nw_k\sum_{j=1}^ka_{k-j}^{(k)}w_j\ge \sum_{k=1}^n\braB{q_{n-k}^{(n)}+\sum_{j=1}^ka_{k-j}^{(k)}}w_k^2>0
\quad \text{for $n\ge1$ if $w_k\not\equiv0$.}
\end{align*}
\end{lemma}
\begin{proof}The definition \eqref{eq: DCO kernels for L1r} implies
$$q_{k-j-1}^{(k-1)}-q_{k-j}^{(k)}=\int_{t_{j-1}}^{t_j}\kbra{\omega_{\alpha}(t_{k-1}-s)-\omega_{\alpha}(t_k-s)}\zd{s}>0,\quad k\ge2,$$
and
\begin{align}\label{eq: sum of kernels}
\sum_{j=1}^ka_{k-j}^{(k)}=\sum_{j=1}^kq_{k-j}^{(k)}-\sum_{j=1}^{k-1}q_{k-j-1}^{(k-1)}=\int_{t_{k-1}}^{t_k}\omega_{\alpha}(s)\zd{s}>0\quad \text{for $k\ge1$.}
\end{align}
Thus we apply the definition \eqref{eq: L1r discrete kernels} to derive that
\begin{align*}
2w_k\sum_{j=1}^ka_{k-j}^{(k)}w_j=&2q_{0}^{(k)}w_k^2-2\sum_{j=1}^{k-1}\brab{q_{k-j-1}^{(k-1)}-q_{k-j}^{(k)}}w_kw_j\\
\ge&2q_{0}^{(k)}w_k^2-\sum_{j=1}^{k-1}\brab{q_{k-j-1}^{(k-1)}-q_{k-j}^{(k)}}\bra{w_k^2+w_j^2}\\
=&\,w_k^2\sum_{j=1}^ka_{k-j}^{(k)}+\sum_{j=1}^kq_{k-j}^{(k)}w_j^2-\sum_{j=1}^{k-1}q_{k-j-1}^{(k-1)}w_j^2\quad \text{for $k\ge1$.}
\end{align*}
This completes the proof.
\end{proof}

\subsection{A Crank-Nicolson type scheme}

We are now ready to propose our numerical scheme. By setting $v:=-\vparl{E}{u}$, one can write the problem \eqref{Cont: AllenCahn RL} into a couple of system
\begin{align}
&\partial_tu=\,{}^R\!\partial_t^{1-\alpha}v, \label{eq: couple problem1}\\
&v=\,\varepsilon^2\Delta u-f(u).\label{eq: couple problem2}
\end{align}
We consider a finite difference approximation in physical domain.
For a positive integer $M_1$, we set the spatial length as $h:=L/M_1$ so that $\bar{\Omega}_{h}:=\big\{\mathbf{x}_{h}=(ih,jh)\,|\,0\leq i,j\leq M_1\}.$
For any grid function $\{v_h\,|\,\mathbf{x}_{h}\in\bar{\Omega}_{h}\}$, we denote
$$\mathbb{V}_{h}:=\big\{v\,|\,v=(v_{j})^{T}\;\;\text{for}\;\;1\leq{j}\leq{M_1},
\;\text{with}\;v_{j}=(v_{i,j})^{T}\;\text{for}\;1\leq{i}\leq{M_1}\big\},$$
where $v^{T}$ is the transpose of the vector $v$.
The maximum norm
$\|v\|_{\infty}:=\max_{\mathbf{x}_{h}\in\bar{\Omega}_{h}}|v_{h}|.$
Let $M:=M_1^2,$ we denote by $D_h$ the $M\times M$ matrix of Laplace operator $\Delta$ subject to periodic boundary conditions.

Now, by applying the L1$_R$ approximation \eqref{eq: L1r formula} in the time domain and a second-order approximation
for the nonlinear term ( see Appendix \ref{append: nonlinear approximation} for details), we obtain
a Crank-Nicolson type scheme in the vector form:
\begin{align}
&\partial_{\tau}u^{n-\frac{1}{2}}=\,\brab{{}^R\!\partial_{\tau}^{1-\alpha}v}^{n-\frac{1}{2}}\quad \text{for $1\leq n\leq N$},\label{eq: L1rCN scheme1}\\
& v^{n-\frac12}=\,\varepsilon^2D_h u^{n-\frac12}-H(u^n,u^{n-1})\quad\text{for $1\leq n\leq N$},\label{eq: L1rCN scheme2}
\end{align}
where the vector $H(u^n,u^{n-1})$ is defined in the element-wise with the Hadamard product ``$\circ$'',
\begin{align}\label{fun: nonlinear approximation}
H(u^n,u^{n-1}):=\frac{1}{3}(u^{n})^{.3}+\frac{1}{2}(u^{n-1})^{.2}\circ u^{n}+\frac{1}{6}(u^{n-1})^{.3}-\frac12\brab{u^{n}+u^{n-1}}.
\end{align}
The constructing procedure of $H(u^n,u^{n-1})$ and its properties are presented in Appendix \ref{append: nonlinear approximation}, and the associated properties of the constructing procedure will be useful for our analysis in later sections.

We close this section by listing some simple properties of the matrix $D_h$ in the following lemma (whose proof is similar as in \cite{HouTangYang:2017Numerical}).
\begin{lemma}\label{lem: matrixD Negative-Condition}
The discrete matrix $D_h$ admits the following properties
\begin{itemize}
  \item [(a)] The discrete matrix $D_h$ is symmetric.
  \item [(b)] For any nonzero $v\in{\mathbb{V}_{h}}$, $v^{T}D_hv\leq{0}$, i.e., the matrix $D_h$ is negative semi-definite.
  \item [(c)] The elements of $D_h=(d_{ij})$ satisfy $d_{ii}=-\max_{i}\sum_{j\neq{i}}|d_{ij}|$ for each $i$.
\end{itemize}
\end{lemma}

\subsection{Discrete variational energy dissipation law}

In this section, we shall establish the discrete variational energy dissipation law for our numerical scheme. To this end, we first define the discrete version of the variational energy. Consider the midpoint rule of the fractional Riemann-Liouville integral operator $\mathcal{I}_t^{\alpha}$ defined by \eqref{Cont: RL integral def},
\begin{align}\label{eq: RL integral formula}
(\mathcal{I}_{t}^{\alpha}v)(t_n)\approx&\, \sum_{k=1}^{n}\int_{t_{k-1}}^{t_k}\omega_{\alpha}(t_n-s)(\Pi_{0,k}v)(s)\zd{s}
=\sum_{k=1}^{n}q_{n-k}^{(n)}v^{k-\frac{1}2}\triangleq(\mathcal{I}_{\tau}^{\alpha}v)^{n}
\end{align}
for $n\ge 0$. Namely, the auxiliary kernels $q_{n-k}^{(n)}$ in \eqref{eq: DCO kernels for L1r}
define a numerical fractional integral $(\mathcal{I}_{\tau}^{\alpha}v)^{n}$.
Notice that the L1$_R$ formula \eqref{eq: L1r formula} yields an alternative formula
for \eqref{Cont: Riemann-Liouville def}, i.e.,
\begin{align}\label{eq: alternative L1r formula}
({}^R\!\partial_{\tau}^{1-\alpha}v)^{n-\frac{1}{2}}=\partial_{\tau}(\mathcal{I}_{\tau}^{\alpha}v)^{n-\frac12}
:=\frac1{\tau_n}\kbra{(\mathcal{I}_{\tau}^{\alpha}v)^n-(\mathcal{I}_{\tau}^{\alpha}v)^{n-1}}\quad\text{for $n\ge 1$}.
\end{align}
We now define the discrete version of our variational energy:
\begin{align*}
\mathcal{E}_{\alpha}[u^n]:=&\,E[u^n]+\frac12h^2\sum_{i,j=1}^{M_1}(\mathcal{I}_{\tau}^{\alpha}v_{ij}^2)^{n}
=E[u^n]+\frac12h^2\sum_{i,j=1}^{M_1}\sum_{k=1}^{n}q_{n-k}^{(n)}\brab{v_{ij}^{k-\frac{1}2}}^2,
\end{align*}
where $v^{k-\frac12}$ represents a numerical approximation at $t_{k-\frac12}$ of the energy variation $\vparl{E}{u},$ and $E[u^n]$ is the discrete counterpart of the free energy \eqref{Cont: free energy}
\begin{align*}
E[u^n]:=&\,h^2\sum_{i,j=1}^{M_1}F(u_{ij}^n)-\frac{1}2\varepsilon^{2}h^2\brab{u^{n}}^TD_{h}u^{n}\quad \text{for $n\ge0$.}
\end{align*}

We are now ready to present the discrete variational energy dissipation law
for our Crank-Nicolson type scheme \eqref{eq: L1rCN scheme1}-\eqref{eq: L1rCN scheme2}.

\begin{theorem}\label{thm: discrete energy dissipation}
The Crank-Nicolson scheme \eqref{eq: L1rCN scheme1}-\eqref{eq: L1rCN scheme2}
admits the variational energy dissipation law unconditionally at the discrete levels, i.e.,
\begin{align*}
\partial_{\tau}\brab{\mathcal{E}_{\alpha}[u]}^{n-\frac12}
+\frac{1}{2\tau_n}\int_{t_{n-1}}^{t_n}\omega_{\alpha}(s)\zd{s}\sum_{i,j=1}^{M_1}h^2\brab{v_{ij}^{n-\frac12}}^2\le0\quad \text{for $n\ge1$.}
\end{align*}
\end{theorem}
\begin{proof}
Taking the $L^2$ inner products of \eqref{eq: L1rCN scheme1}-\eqref{eq: L1rCN scheme2} with $\tau_n\brat{v^{n-\frac12}}^T$
and $-\brat{\diff u^n}^T$, respectively, and adding up the two resulting equalities, we obtain
\begin{align}\label{thmproof:discrete energy}
\sum_{i,j=1}^{M_1}h^2v_{ij}^{n-\frac12}\sum_{k=1}^na_{n-k}^{(n)}v_{ij}^{k-\frac12}
-&\,\varepsilon^{2}h^2\brab{\diff u^n}^TD_{h}u^{n-\frac12}\nonumber\\
+&\,\sum_{i,j=1}^{M_1}h^2H(u_{ij}^n,u_{ij}^{n-1})\brab{\diff u_{ij}^n}=0\quad \text{for $n\ge1$.}
\end{align}
With the help of  Lemma \ref{lem: matrixD Negative-Condition} (a)-(b), it is easy to show that
\begin{align}
\brab{\diff u^n}^TD_{h}u^{n-\frac12}=&\,
\frac{1}{2}\brab{\diff u^n}^TD_{h}\brat{u^{n}+u^{n-1}} \nonumber\\
=&\,\frac{1}2\brab{u^{n}}^TD_{h}u^{n}-\frac{1}2\brab{u^{n-1}}^TD_{h}u^{n-1}
\quad\text{for $n\ge1$.}
\end{align}
Moreover, by taking $a=u_{ij}^n$ and $b=u_{ij}^{n-1}$ in the equality \eqref{fun: nonlinear coercivity} one has
\begin{align}
H(u_{ij}^n,u_{ij}^{n-1})\brab{\diff u_{ij}^n}\ge F(u_{ij}^n)-F(u_{ij}^{n-1})\quad\text{for $n\ge1$.}
\end{align}
By using Lemma \ref{lem: L1r positive definite} and the formulas \eqref{eq: RL integral formula}-\eqref{eq: alternative L1r formula}, one has
\begin{align}\label{eq11}
v_{ij}^{n-\frac12}\sum_{k=1}^na_{n-k}^{(n)}v_{ij}^{k-\frac12}\ge&\,
\frac12(\mathcal{I}_{\tau}^{\alpha}v_{ij}^2)^{n}
-\frac12(\mathcal{I}_{\tau}^{\alpha}v_{ij}^2)^{n-1}+\frac12\brab{v_{ij}^{n-\frac12}}^2\sum_{k=1}^na_{n-k}^{(n)}\nonumber\\
=&\,
\frac12(\mathcal{I}_{\tau}^{\alpha}v_{ij}^2)^{n}
-\frac12(\mathcal{I}_{\tau}^{\alpha}v_{ij}^2)^{n-1}+\frac{1}2\int_{t_{n-1}}^{t_n}\omega_{\alpha}(s)\zd{s}\cdot\brab{v_{ij}^{n-\frac12}}^2,
\end{align}
where \eqref{eq: sum of kernels} was used.  Thus the claimed result follows from \eqref{thmproof:discrete energy}-\eqref{eq11} immediately.
\end{proof}

Notice that as the fractional order $\alpha\rightarrow1$, the definition \eqref{eq: DCO kernels for L1r} yields
$q_{n-k}^{(n)}\rightarrow\tau_k.$ Moreover, the numerical fractional integral
and the L1$_R$ formula yield
$$(\mathcal{I}_{\tau}^{\alpha}v)^{n}\rightarrow \sum_{k=1}^{n}\tau_kv^{k-\frac{1}2}, \quad \brab{{}^R\!\partial_{\tau}^{1-\alpha}v}^{n-\frac{1}{2}}\rightarrow v^{n-\frac{1}{2}}.$$
Consequently, the discrete variational energy dissipation law in Theorem \ref{thm: discrete energy dissipation} becomes
\begin{align*}
\partial_{\tau}\brab{E[u]}^{n-\frac12}+\sum_{i,j=1}^{M_1}h^2\brab{v_{ij}^{n-\frac12}}^2\le 0\quad \text{as $\alpha\rightarrow1$.}
\end{align*}
This recovers the standard discrete energy dissipation law of the classical Allen-Cahn equation. Thus, the discrete variational energy stability (Theorem 2.1) is asymptotically preserving in the $\alpha \rightarrow 1$ limit.

\section{Unique solvability and discrete maximum bound principle}
\setcounter{equation}{0}
This section will be devoted to the unique solvability and discrete maximum bound principle of our scheme.
To this end, we shall first introduce some analysis tools  including the
discrete orthogonal convolution (DOC) kernels and discrete complementary-to-orthogonal (DCO) kernels.
\subsection{DOC and DCO kernels}
We first introduce a class of discrete orthogonal convolution (DOC) kernels $\theta_{n-j}^{(n)}$ via the following discrete orthogonal identity (with respect to the discrete kernels $a_{n-k}^{(n)}$ in (\ref{eq: L1r discrete kernels}))
\begin{align}\label{eq: orthogonal identity}
\sum_{j=k}^{n}\theta_{n-j}^{(n)}a^{(j)}_{j-k}\equiv\delta_{nk}
\quad\text{for $\forall\;1\leq k\le n$,}
\end{align}
where $\delta_{nk}$ is the Kronecker delta symbol.
Notice that the DOC kernels can be defined via a recursive procedure
\begin{align}\label{eq: orthogonal procedure}
\theta_{0}^{(n)}:=\frac{1}{a^{(n)}_{0}}\quad\text{and}\quad
\theta_{n-k}^{(n)}:=-\frac1{a^{(k)}_{0}}\sum_{j=k+1}^{n}\theta_{n-j}^{(n)}a^{(j)}_{j-k}
\quad\text{for $1\leq k\le n-1$.}
\end{align}
This type of discrete kernels has been used in \cite{LiaoZhang:2020} for analyzing the nonuniform BDF2 scheme for linear diffusion problems.
Here we consider the DOC kernels $\theta_{n-j}^{(n)}$ of the L1$_R$ discrete kernels \eqref{eq: L1r discrete kernels}
that satisfy
 $$a^{(n)}_{0}>0\quad\text{and}\quad a^{(n)}_{j}<0\quad\text{for $1\leq j\le n-1$.}$$

We shall show that, by using the DOC kernels  $\theta_{n-j}^{(n)},$
the Crank-Nicolson scheme \eqref{eq: L1rCN scheme1}-\eqref{eq: L1rCN scheme2}
in the Riemann-Liouville form can be reformulated into an equivalent form in the Caputo form, see \eqref{eq: L1r caputo scheme} below.
Then, the unique solvability and discrete maximum principle of
Crank-Nicolson scheme can be performed via the equivalent form \eqref{eq: L1r caputo scheme}.

Furthermore, the original discrete form \eqref{eq: L1rCN scheme1}-\eqref{eq: L1rCN scheme2} can also be recovered from \eqref{eq: L1r caputo scheme}
by using the L1$_R$ discrete kernels $a_{n-j}^{(n)}.$ This seems to be the first discrete transformation between two discrete fractional derivatives, and  the discrete duality between the transform and the inverse transform relies on the following mutual orthogonality.

\begin{lemma}\cite[Lemma 2.1]{LiaoTangZhou:2020doc}\label{lem: Mutual orthogonality}
The discrete convolution kernels $a^{(n)}_{n-j}$ and the corresponding DOC kernels $\theta_{n-j}^{(n)}$ are mutually orthogonal, that is,
\begin{align}\label{eq: Mutual orthogonal identity}
\sum_{j=k}^na^{(n)}_{n-j}\theta_{j-k}^{(j)}\equiv\delta_{nk}\quad\text{and}\quad \sum_{j=k}^{n}\theta_{n-j}^{(n)}a^{(j)}_{j-k}\equiv\delta_{nk}
\quad\text{for $1\le k\le n$.}
\end{align}
\end{lemma}
We now list some useful properties of the DOC kernels $\theta_{n-j}^{(n)}.$
\begin{lemma}\label{lem: lower bound A0}
For fixed $n\ge1$, the DOC kernels $\theta_{n-j}^{(n)}$ are positive and
satisfy
\begin{align}\label{lower bounds}
\theta_{0}^{(n)}=\frac{1}{\omega_{1+\alpha}(\tau_{n})}\quad\text{and}\quad
\theta_{0}^{(n)}-\theta_{1}^{(n)}>\frac{\omega_{\alpha}(r_n+1)}{\omega_{1+\alpha}(\tau_{n})\omega_{1+\alpha}(1)}\,.
\end{align}
\end{lemma}
\begin{proof}
By the definitions \eqref{eq: L1r discrete kernels}-\eqref{eq: DCO kernels for L1r}, we have $a^{(n)}_{0}=q_{0}^{(n)}=\omega_{1+\alpha}(\tau_n)$
so that the procedure \eqref{eq: orthogonal procedure} yields
$$\theta^{(n)}_{0}=\frac1{a_{0}^{(n)}}=\frac{1}{\omega_{1+\alpha}(\tau_{n})}>0\quad\text{for $n\ge1$.}$$
The positivity of DOC kernels $\theta^{(n)}_{j}$ can be verified by a simple induction argument.
Assume that $\theta^{(n)}_{j}>0$ for $0\le j\le m-1$ $(m\ge1)$.
By the definition \eqref{eq: L1r discrete kernels}, one has 
\begin{align*}
a^{(j)}_{j-k}=\int_{t_{j-1}}^{t_{j}}\zd{t}\int_{t_{k-1}}^{t_{k}}\omega_{\alpha-1}(t-s)\zd{s}<0
\quad\text{for $j\ge k+1\ge2$.}
\end{align*}
Then the recursive procedure \eqref{eq: orthogonal procedure} gives
\begin{align*}
\theta_{m}^{(n)}=-\frac{1}{a^{(n-m)}_{0}}\sum_{\ell=0}^{m-1}\theta_{\ell}^{(n)}a^{(n-\ell)}_{m-\ell}>0
\quad\text{for $1\leq m\le n-1$.}
\end{align*}
This confirms that $\theta^{(n)}_{j}>0$ for $0\le j\le n-1.$

Moreover, we have $\theta_{1}^{(n)}=-\theta_{0}^{(n)}a^{(n)}_{1}/a^{(n-1)}_{0}$
so that the definition \eqref{eq: DCO kernels for L1r} yields
\begin{align*}
\theta_{0}^{(n)}-\theta_{1}^{(n)}=&\,\theta^{(n)}_{0}+\theta_{0}^{(n)}\frac{a^{(n)}_{1}}{a^{(n-1)}_{0}}
=\frac{q^{(n)}_{1}}{q^{(n)}_{0}q^{(n-1)}_{0}}\\
=&\,\frac{\omega_{1+\alpha}(\tau_n+\tau_{n-1})-\omega_{1+\alpha}(\tau_{n})}
{\omega_{1+\alpha}(\tau_{n})\omega_{1+\alpha}(\tau_{n-1})}
=\frac{\omega_{1+\alpha}(r_n+1)-\omega_{1+\alpha}(r_{n})}
{\omega_{1+\alpha}(\tau_{n})\omega_{1+\alpha}(1)}.
\end{align*}
This yields the claimed lower bound (\ref{lower bounds}), and the proof is completed.
\end{proof}

\begin{lemma}\label{lem: DOC positive decreasing}
For any fixed $n\ge2$, the DOC kernels $\theta_{n-j}^{(n)}$ are monotonously decreasing, that is,  $$\theta_{0}^{(n)}>\theta_{1}^{(n)}>\cdots>\theta_{n-1}^{(n)}>0.$$
\end{lemma}
\begin{proof}
Applying the definitions \eqref{eq: L1r discrete kernels} and \eqref{eq: orthogonal procedure}, one has $\theta_{0}^{(n)}q^{(n)}_{0}=1$ and
\begin{align*}
\theta_{n-k}^{(n)}q^{(k)}_{0}
=&\,-\sum_{j=k+1}^{n}\theta_{n-j}^{(n)}\brab{q^{(j)}_{j-k}-q^{(j-1)}_{j-k-1}}\\
=&\,-\sum_{j=k+1}^{n}\theta_{n-j}^{(n)}q^{(j)}_{j-k}+\sum_{j=k}^{n-1}\theta_{n-j-1}^{(n)}q^{(j)}_{j-k}
\quad\text{for $1\leq k\le n-1$,}
\end{align*}
or
\begin{align*}
\theta_{0}^{(n)}q^{(n)}_{n-k}
=&\,-\sum_{j=k}^{n-1}\brab{\theta_{n-j}^{(n)}-\theta_{n-j-1}^{(n)}}q^{(j)}_{j-k}
\quad\text{for $1\leq k\le n-1$.}
\end{align*}
Consider an auxiliary class of discrete kernels $\zeta^{(n)}_{n-j}$ defined by
\begin{align*}
\zeta^{(n)}_{0}:=\theta_{0}^{(n)}\quad\text{and}\quad
\zeta^{(n)}_{n-j}:=\theta_{n-j}^{(n)}-\theta_{n-j-1}^{(n)}\quad\text{for $1\leq j\le n-1$.}
\end{align*}
Then it is easy to find that
\begin{align*}
\sum_{j=k}^{n}\zeta_{n-j}^{(n)}q^{(j)}_{j-k}=\delta_{nk}\quad\text{for $1\leq k\le n$,}
\end{align*}
that is, the kernels $\zeta^{(n)}_{n-j}$ are orthogonal to $q^{(n)}_{n-k}=\int_{t_{k-1}}^{t_k}\omega_{\alpha}(t_n-s)\zd{s}$.
By following the proof of \cite[Proposition 4.1]{LiaoTangZhou:2020doc},  it is easy to check that
\begin{align*}
q_{j}^{(n)}>0,\quad q^{(n-1)}_{j-1}>q^{(n)}_{j}\quad\text{and}\quad q^{(n-1)}_{j-1}q^{(n)}_{j+1}>q^{(n-1)}_{j}q^{(n)}_{j}.
\end{align*}
Then \cite[Lemma 2.3]{LiaoTangZhou:2020doc} implies that the corresponding orthogonal kernels $\zeta^{(n)}_{n-j}$ satisfy
\begin{align*}
\zeta^{(n)}_{n-j}=\theta_{n-j}^{(n)}-\theta_{n-j-1}^{(n)}<0\quad\text{for $1\leq j\le n-1$}\quad
\text{and}\quad  \theta_{n-1}^{(n)}=\sum_{j=1}^n\zeta^{(n)}_{n-j}>0.
\end{align*}
They implies the claimed property and complete the proof.
\end{proof}

\begin{lemma}\label{lem: complementary DCO}
The discrete convolution kernels $q_{n-j}^{(n)}$ in \eqref{eq: DCO kernels for L1r} are
complementary to the DOC kernels $\theta_{n-j}^{(n)}$ in \eqref{eq: orthogonal procedure} in the sense that
\begin{align*}
&\sum_{j=k}^{n}q_{n-j}^{(n)}\theta_{j-k}^{(j)}\equiv 1\quad \text{for $ 1\leq k\le n$.}
\end{align*}
\end{lemma}
\begin{proof}Inserting the definition \eqref{eq: L1r discrete kernels} into the first identity of \eqref{eq: Mutual orthogonal identity}
arrives
\begin{align}\label{eq: complementary procedure}
q^{(n)}_{0}\theta_{n-k}^{(n)}+\sum_{j=k}^{n-1}\brab{q_{n-j}^{(n)}-q_{n-j-1}^{(n-1)}}\theta_{j-k}^{(j)}
\equiv\delta_{nk}\quad \text{for $ 1\leq k\le n$,}
\end{align}
which implies
\begin{align*}
\sum_{j=k}^{n}q_{n-j}^{(n)}\theta_{j-k}^{(j)}\equiv\sum_{j=k}^{n-1}q_{n-j-1}^{(n-1)}\theta_{j-k}^{(j)}+\delta_{nk}
\quad \text{for $ 1\leq k\le n$.}
\end{align*}
Let $\Xi_k^{(n)}:=\sum_{j=k}^{n}q_{n-j}^{(n)}\theta_{j-k}^{(j)}$ for $1\le k\le n$. One has
\begin{align*}
\Xi_n^{(n)}=1\quad\text{and}\quad\Xi_k^{(n)}=\Xi_k^{(n-1)}\quad \text{for $ 1\leq k\le n-1$.}
\end{align*}
A simple induction yields $\Xi_k^{(n)}\equiv1$ for $1\le k\le n$  and completes the proof.
\end{proof}

\begin{figure}[!ht]
\caption{The relationship diagram between different classes of discrete kernels.}
  \centering
  \includegraphics[width=4in]{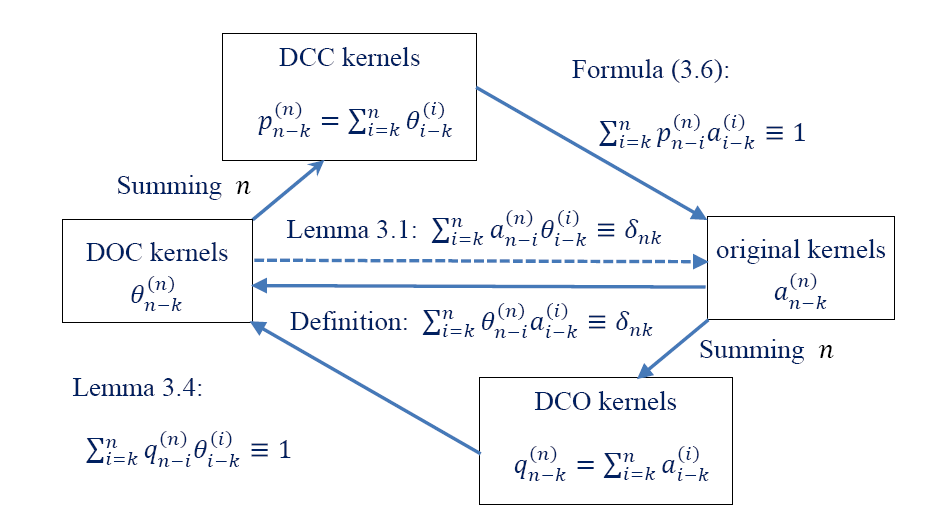}
  \label{fig: DOC, DCO and DCC relation}
\end{figure}

Since the kernels $q_{n-j}^{(n)}$ are  complementary to the DOC kernels $\theta_{n-j}^{(n)},$
we call $q_{n-j}^{(n)}$ as the discrete complementary-to-orthogonal (DCO) kernels. This terminology is used to distinguish it from
  the discrete complementary convolution (DCC) kernels $p_{n-j}^{(n)}$
  which are complementary to the original discrete kernels $a_{n-j}^{(n)}$, i.e.,
\begin{align}\label{eq: complementary DPA}
\sum_{j=k}^{n}p_{n-j}^{(n)}a_{j-k}^{(j)}\equiv 1\quad \text{for $ 1\leq k\le n$.}
\end{align}
We present in Figure \ref{fig: DOC, DCO and DCC relation} the relationships between
the mentioned discrete convolution kernels.
We notice that the DCC kernels $p_{n-j}^{(n)}$ were originally introduced in \cite{LiaoLiZhang:2018,LiaoMcLeanZhang:2019}
for analyzing the direct approximations of Caputo derivative.
Here, we shall use the newly introduced DCO kernels $q_{n-j}^{(n)}$ to analyze the direct approximations of the Riemann-Liouville derivative \eqref{Cont: Riemann-Liouville def}.

\subsection{An equivalent formula and unique solvability}

We now derive an equivalent formula for our scheme \eqref{eq: L1rCN scheme1}-\eqref{eq: L1rCN scheme2}.
By using the definition \eqref{eq: L1r formula} of L1$_R$ formula, we can write the equation \eqref{eq: L1rCN scheme1} as
\begin{align*}
\diff u^{j}=\tau_j({}^R\!\partial_{\tau}^{1-\alpha}v)^{j-\frac{1}{2}}
=\sum_{k=1}^{j}a_{j-k}^{(j)}v^{k-\frac{1}2}\quad\text{for $1\leq j\leq N$.}
\end{align*}
Multiplying both sides of the above equation by the DOC kernels
$\theta_{n-j}^{(n)}$, and summing $j$ from $j=1$ to $n$, we obtain
\begin{align*}
\sum_{j=1}^n\theta_{n-j}^{(n)}\diff u^{j}=&\,
\sum_{j=1}^n\theta_{n-j}^{(n)}\sum_{k=1}^{j}a_{j-k}^{(j)}v^{k-\frac{1}2}
=\sum_{k=1}^{n}v^{k-\frac{1}2}
\sum_{j=k}^n\theta_{n-j}^{(n)}a_{j-k}^{(j)}\nonumber\\
=&\,\sum_{k=1}^{n}v^{k-\frac{1}2}\delta_{nk}=v^{n-\frac{1}2}
\quad\text{for $1\leq n\leq N$,}
\end{align*}
where the summation order was exchanged in the second equality
and the discrete orthogonal identity \eqref{eq: orthogonal identity} was used in the third equality.
Then the equation \eqref{eq: L1rCN scheme2} gives an equivalent form of the Crank-Nicolson scheme
\begin{align}\label{eq: L1r caputo scheme}
\sum_{j=1}^n\theta_{n-j}^{(n)}\diff u^{j}=\varepsilon^2D_hu^{n-\frac{1}{2}}-H(u^{n},u^{n-1})\quad\text{for $1\leq n\leq N$},
\end{align}
where the vector $H(u^n,u^{n-1})$ is defined by \eqref{fun: nonlinear approximation}.

The formulation \eqref{eq: L1r caputo scheme} looks like a direct approximation of the original equation \eqref{Cont: TFAC}
by approximating the Caputo derivative $\partial_t^{\alpha}u$ with
\begin{align}\label{eq: L1r caputo formula}
(\partial_t^{\alpha}u)(t_{n-\frac12})\approx\sum_{j=1}^n\theta_{n-j}^{(n)}\diff u^{j}.
\end{align}
In this sense, the DOC kernels $\theta_{n-j}^{(n)}$ define a ``new" discrete Caputo derivative.
According to Lemma \ref{lem: DOC positive decreasing}, the corresponding discrete kernels $\theta_{n-j}^{(n)}$
are positive and monotonously decreasing on nonuniform time meshes, as the case of L1 formula \cite{LiaoLiZhang:2018,LiaoMcLeanZhang:2019,LiaoYanZhang:2019}.
Nonetheless, we view the formula \eqref{eq: L1r caputo formula} as an indirect approximation
that admits different approximation accuracy compared to the original L1 formula
with the approximation error of $2-\alpha$.

Next, we shall proof the solvability and discrete maximum bound principle via the new form \eqref{eq: L1r caputo scheme}. To this end, we shall also need the following lemma for which the proof is similar as in \cite[Lemma 3.2]{HouTangYang:2017Numerical}.
\begin{lemma}\label{lem: Matrix-Inf-Norm}
Let the elements of a real matrix $B=(b_{ij})_{M\times{M}}$ fulfill
$b_{ii}=-\max_{i}\sum_{j\neq{i}}|b_{ij}|\,.$
For any parameters $a,c>0$ and $U,V\in{\mathbb{R}^{M}}$, it holds that
$$\mynormb{\bra{aI-B}V}_{\infty}\geq{a}\mynormb{V}_{\infty}$$
and
\begin{align*}
\mynormb{\bra{aI-B}V+U^{.2}\circ V+cV^{.3}}_{\infty}\geq{a}\mynormb{V}_{\infty}+\mynormb{U}_{\infty}^{2}\mynormb{V}_{\infty}+c\mynormb{V}_{\infty}^{3}.
\end{align*}
\end{lemma}
Now we are ready to present the unique solvability of our scheme.
\begin{theorem}\label{thm: unique solvability}
The nonlinear Crank-Nicolson scheme \eqref{eq: L1r caputo scheme} or \eqref{eq: L1rCN scheme1}-\eqref{eq: L1rCN scheme2} is uniquely solvable
if  the time-step size satisfies $$\tau< \sqrt[\alpha]{2\Gamma(1+\alpha)}.$$
\end{theorem}
\begin{proof}
We rewrite the nonlinear scheme \eqref{eq: L1r caputo scheme} into
\begin{align*}
G_hu^{n}+\frac13(u^{n})^{.3}=\mathrm{G}_0(u^{n-1}),\quad
{n}\geq{1},
\end{align*}
where
$$G_h:=\braB{\theta_0^{(n)}-\frac12+\frac{1}{2}(u^{n-1})^{.2}}I-\frac{\varepsilon^{2}}2D_{h}$$
and
\begin{align*}
\mathrm{G}_0(u^{n-1}):=&\,
\frac{1}2\brat{I+\varepsilon^{2}D_{h}}u^{n-1}-\frac16(u^{n-1})^{.3}
+\sum_{k=1}^{n-1}\brab{\theta_{n-k-1}^{(n)}-\theta_{n-k}^{(n)}}u^k+\theta_{n-1}^{(n)}u^0
\quad\text{for $n\geq1$}.
\end{align*}
If the maximum step size $\tau<\sqrt[\alpha]{2\Gamma(1+\alpha)}$, Lemma \ref{lem: lower bound A0}
shows that $$\theta_0^{(n)}=\Gamma(1+\alpha)\tau_n^{-\alpha}> 1/2.$$
Then by Lemma \ref{lem: matrixD Negative-Condition} (b), the symmetric matrix $G_h$ is positive definite.
Thus, the solution of nonlinear equations solves
  \begin{align*}
u^n=\arg \min_{w\in \mathbb{V}_{h}}\left\{\frac12w^TG_hw+\frac{1}{12}\sum_{k=1}^{M}w_k^{4}-w^T\mathrm{G}_0(u^{n-1})\right\}\quad
\text{for $n\geq1$.}
\end{align*}
The strict convexity of objective function implies the unique solvability of \eqref{eq: L1r caputo scheme}.
\end{proof}

\subsection{Discrete maximum bound principle}
We next show that our scheme admits the discrete maximum bound principle.
\begin{theorem}\label{thm:Dis-Max-Principle}
Assume that the time-step size satisfies
\begin{align}\label{Constraint: Dis_Max Time-Steps}
\tau_n\le \sqrt[\alpha]{\min\Big\{\frac12,\frac{h^2}{2\varepsilon^{2}}\Big\}\frac{\alpha\,\Gamma(1+\alpha)}{(1+r_n)^{1-\alpha}}}\,.
\end{align}
Then, the Crank-Nicolson scheme \eqref{eq: L1r caputo scheme} or \eqref{eq: L1rCN scheme1}-\eqref{eq: L1rCN scheme2}
preserves the maximum bound principle at the discrete levels,
that is,
\begin{align*}
\mynormb{u^{k}}_{\infty}\leq{1}\;\;\text{for $1\leq{k}\leq{N}$}\quad\text{if}\;\;\mynormb{u^{0}}_{\infty}\leq{1}.
\end{align*}
\end{theorem}
\begin{proof}According to Theorem \ref{thm: unique solvability},
the time-step restriction \eqref{Constraint: Dis_Max Time-Steps} ensures the solvability of \eqref{eq: L1r caputo scheme}
since $(1+r_n)^{\alpha-1}<1$.
Now we consider a mathematical induction proof.
Obviously, the claimed inequality holds for $n=0$.
For $1\leq n\le N$, assume that
\begin{align}\label{thmproof: Max-Principle inductionAssume}
\mynormb{u^{k}}_{\infty}\leq{1}\quad\text{for $0\leq{k}\leq{n-1}.$}
\end{align}
It remains to verify that $\mynormb{u^{n}}_{\infty}\leq{1}$.
Note that
\begin{align*}
\sum_{j=1}^n\theta_{n-j}^{(n)}\diff u^{j}
=\theta_{0}^{(n)}u^{n}
-\brab{\theta_{0}^{(n)}-\theta_{1}^{(n)}}u^{n-1}-\mathcal{L}^{n-2}(u),
\end{align*}
where $\mathcal{L}^{n-2}(u)$ is given by
\begin{align}\label{thmproof: Max-Principle histroy operator}
\mathcal{L}^{n-2}(u):=\sum_{k=1}^{n-2}\big(\theta_{n-k-1}^{(n)}-\theta_{n-k}^{(n)}\big)u^{k}+\theta_{n-1}^{(n)}u^{0}.
\end{align}
Then the scheme \eqref{eq: L1r caputo scheme} can be formulated as follows
\begin{align}\label{thmproof: Max-Principle Scheme-Variant}
&\,\braB{\theta^{(n)}_{0}-\frac12+\frac{1}{2}(u^{n-1})^{.2}-\frac{\varepsilon^{2}}2D_{h}}u^{n}
+\frac13(u^{n})^{.3}\nonumber\\
&\,\hspace{2cm}=\brab{\theta_{0}^{(n)}-\theta_{1}^{(n)}}u^{n-1}+\frac{\varepsilon^{2}}2D_{h}u^{n-1}
+\frac12u^{n-1}-\frac16(u^{n-1})^{.3}+\mathcal{L}^{n-2}(u)\nonumber\\
&\,\hspace{2cm}=\mathbf{M}_hu^{n-1}+\frac16\kbra{3u^{n-1}-(u^{n-1})^{.3}}+\mathcal{L}^{n-2}(u),
\end{align}
where the matrix $\mathbf{M}_h$ is defined by
\begin{align}\label{thmproof: Max-Principle matrix Qh}
\mathbf{M}_h:=\brab{\theta_{0}^{(n)}-\theta_{1}^{(n)}}I+\frac{\varepsilon^{2}}2D_{h}.
\end{align}

For the first term of the right hand side of \eqref{thmproof: Max-Principle Scheme-Variant}, it is easy to check that the matrix $\mathbf{M}_h=(m_{ij})$ satisfies $m_{ij}\ge0$ for $i\neq j$,
\begin{align*}
m_{ii}=\theta_{0}^{(n)}-\theta_{1}^{(n)}-\frac{2\varepsilon^{2}}{h^2}
\quad\text{and}\quad
\max_{i}\sum_{j}m_{ij}\le \theta_{0}^{(n)}-\theta_{1}^{(n)}.
\end{align*}
Assuming that $\tau_n\le \sqrt[\alpha]{\frac{h^2}{2\varepsilon^{2}}\Gamma^2(1+\alpha)\omega_{\alpha}(1+r_n)}$,
we apply Lemma \ref{lem: lower bound A0} to find
\begin{align*}
\theta_{0}^{(n)}-\theta_{1}^{(n)}>
\Gamma^2(1+\alpha)\omega_{\alpha}(1+r_n)\tau_n^{-\alpha}\ge\frac{2\varepsilon^{2}}{h^2}
\end{align*}
or $m_{ii}\ge0$.  Thus all elements of $\mathbf{M}_h$ are nonnegative and
\begin{align}\label{thmproof: Max estimate0-matrix Mh}
\mynormb{\mathbf{M}_h}_{\infty}=\max_{i}\sum_{j}\abs{m_{ij}}=\max_{i}\sum_{j}m_{ij}\le \theta_{0}^{(n)}-\theta_{1}^{(n)}.
\end{align}
Consequently, the induction hypothesis \eqref{thmproof: Max-Principle inductionAssume} yields
\begin{align}\label{thmproof: Max estimate1-Scheme-Variant}
\mynormb{\mathbf{M}_hu^{n-1}}_{\infty}\le \mynormb{\mathbf{M}_h}_{\infty}\mynormb{u^{n-1}}_{\infty}
\le \frac12\brab{\theta_{0}^{(n)}-\theta_{1}^{(n)}}\bra{1+\mynormb{u^{n-1}}_{\infty}}.
\end{align}
Since $\abs{3z-z^{3}}\le 2$ for any $z\in[-1,1]$, the induction hypothesis \eqref{thmproof: Max-Principle inductionAssume} yields
\begin{align}\label{thmproof: Max estimate2-Scheme-Variant}
\frac16\mynormb{3u^{n-1}-(u^{n-1})^{.3}}_{\infty}\le \frac13.
\end{align}
For the last term $\mathcal{L}^{n-2}(u)$ in \eqref{thmproof: Max-Principle Scheme-Variant}, the decreasing property in Lemma \ref{lem: DOC positive decreasing}
 and the induction hypothesis \eqref{thmproof: Max-Principle inductionAssume} lead to
\begin{align}\label{thmproof: Max estimate3-Scheme-Variant}
\mynormb{\mathcal{L}^{n-2}(u)}_{\infty}\leq \sum_{k=1}^{n-2}\brab{\theta_{n-k-1}^{(n)}-\theta_{n-k}^{(n)}}\mynormb{u^{k}}_{\infty}
+\theta_{n-1}^{(n)}\mynormb{u^{0}}_{\infty}\leq \theta_{1}^{(n)}.
\end{align}

Moreover, the time-step restriction \eqref{Constraint: Dis_Max Time-Steps} implies $\tau_n< \sqrt[\alpha]{\Gamma(1+\alpha)/2}$
and Lemma \ref{lem: lower bound A0} gives $\theta^{(n)}_{0}>2$.
Then by using Lemmas \ref{lem: matrixD Negative-Condition} and \ref{lem: Matrix-Inf-Norm}, one can
bound the left hand side of \eqref{thmproof: Max-Principle Scheme-Variant} by
\begin{align*}
&\mynormB{\brab{\theta^{(n)}_{0}-\frac12-\frac{\varepsilon^{2}}{2}D_{h}}u^{n}+\frac{1}{2}(u^{n-1})^{.2}\circ u^{n}
+\frac13(u^{n})^{.3}}_{\infty}\\
&\hspace{3.7cm}\ge \brab{\theta^{(n)}_{0}-\frac12}\mynormb{u^{n}}_{\infty}+\frac12\mynormb{u^{n-1}}_{\infty}^2\mynormb{u^{n}}_{\infty}
+\frac13\mynormb{u^{n}}_{\infty}^3.
\end{align*}
Consequently, by collecting the estimates \eqref{thmproof: Max estimate1-Scheme-Variant}--\eqref{thmproof: Max estimate3-Scheme-Variant},
it follows from \eqref{thmproof: Max-Principle Scheme-Variant} that
\begin{align}\label{thmproof: Max estimate4-Scheme-Variant}
\brab{\theta^{(n)}_{0}-\frac12}&\,\mynormb{u^{n}}_{\infty}+\frac12\mynormb{u^{n-1}}_{\infty}^2\mynormb{u^{n}}_{\infty}
+\frac13\mynormb{u^{n}}_{\infty}^3\nonumber\\
\le&\,\mynormb{\mathbf{M}_hu^{n-1}+\frac12\brab{\theta_{0}^{(n)}-\theta_{1}^{(n)}}u^{n-1}+\frac16
\kbrab{3u^{n-1}-(u^{n-1})^{.3}}+\mathcal{L}^{n-2}(u)}_{\infty}\nonumber\\
\le&\,\mynormb{\mathbf{M}_hu^{n-1}}_{\infty}
+\frac16\mynormb{3u^{n-1}-(u^{n-1})^{.3}}_{\infty}
+\mynormb{\mathcal{L}^{n-2}(u)}_{\infty}\nonumber\\
\le&\,\frac12\brab{\theta_{0}^{(n)}-\theta_{1}^{(n)}}\bra{1+\mynormb{u^{n-1}}_{\infty}}+\frac13+\theta_{1}^{(n)}.
\end{align}

Assuming that $\tau_n\le \sqrt[\alpha]{\frac{1}{2}\Gamma^2(1+\alpha)\omega_{\alpha}(1+r_n)}$
such that $\theta_{0}^{(n)}-\theta_{1}^{(n)}>2$, we prove $\mynormb{u^{n}}_{\infty}\le1$ by contradiction.
If $\mynormb{u^{n}}_{\infty}>1$,  the above inequality \eqref{thmproof: Max estimate4-Scheme-Variant} requires
\begin{align*}
\theta^{(n)}_{0}-\frac12+\frac12\mynormb{u^{n-1}}_{\infty}^2+\frac13
<\frac12\brab{\theta_{0}^{(n)}-\theta_{1}^{(n)}}\bra{1+\mynormb{u^{n-1}}_{\infty}}
+\frac13+\theta_{1}^{(n)}
\end{align*}
because the following function
\[
g(z):=\brab{\theta^{(n)}_{0}-\frac12+\frac12\mynormb{u^{n-1}}_{\infty}^2}z
+\frac13z^{3}\quad\text{for $z>0$}
\]
is monotonically increasing. It follows that
\begin{align*}
\frac12\bra{1-\mynormb{u^{n-1}}_{\infty}}^2
<\frac12\brab{\theta_{0}^{(n)}-\theta_{1}^{(n)}-1-\mynormb{u^{n-1}}_{\infty}}\bra{1-\mynormb{u^{n-1}}_{\infty}}<0,
\end{align*}
which yields a contradiction. Thus, the assumption $\mynormb{u^{n}}_{\infty}>1$ is invalid
and the claimed result holds for $k=n$. This completes the proof.
\end{proof}

Note that, the maximum time-step restriction \eqref{Constraint: Dis_Max Time-Steps} is only a sufficient condition to
ensure the discrete maximum principle, see Example \ref{exam:coarsen dynamic}.
In the time-fractional Allen-Cahn equation \eqref{Cont: TFAC},
the coefficient $\varepsilon\ll 1$ represents the width of diffusive interface. Always,
we should choose a small space length $h=O(\varepsilon)$ to track the moving interface.
So, in most situations, the restriction \eqref{Constraint: Dis_Max Time-Steps} is practically reasonable because
it is approximately equivalent to
$$\tau_n\leq\sqrt[\alpha]{\frac{\alpha\,\Gamma(1+\alpha)}{2(1+r_n)^{1-\alpha}}}\rightarrow \frac12\quad\text{as $\alpha\rightarrow1$.}$$
As expected, this restriction \eqref{Constraint: Dis_Max Time-Steps} requires small time steps
for large step ratios $r_n$ or small fractional orders $\alpha$, see similar conditions in \cite{LiaoTangZhou:2020jcp}.
On the other hand, this time-step condition is sharp in the sense that
it is compatible with the previous restriction \cite[Theorem 1]{HouTangYang:2017Numerical} ensuring the discrete maximum principle of
Crank-Nicolson scheme for the classical Allen-Cahn equation.

\section{Numerical experiments}
\setcounter{equation}{0}

In this section, we shall present several numerical examples to support our theoretical findings.
To speed up our numerical computations, we shall use the fast L1$_R$ algorithm described in Appendix \ref{append: fast L1r algorithm},  with
an absolute tolerance error $\epsilon=10^{-12}$ and a cut-off time $\Delta{t}=10^{-12}$.

\subsection{Accuracy verification}
We first show the accuracy of our scheme. Notice that it was shown that the L$1_R$ formula \eqref{eq: L1r formula} has been investigated in \cite{Mustapha:2011,MustaphaAlMutawa:2012}
for linear subdiffusion problems, and the approximation order is shown to be $1+\alpha.$ Thus, we also expect a $(1+\alpha)$-order rate of convergence.
\begin{example}\label{exam:accuracy test}
Consider the exterior-forced model
$$\partial_t u=-{}^{R}\!\partial_t^{1-\alpha}\brab{\frac{\delta E}{\delta u}}+g(\mathbf{x},t)$$
on the space-time domain $(0,1)^{2}\times(0,1]$ with an interfacial coefficient $\varepsilon=0.1$.
We choose a exterior force $g(\mathbf{x},t)$ and a parameter $\sigma\in(0,1)$
such that the model has an exact solution $u=\omega_{1+\sigma}(t)\sin(2\pi x)\sin(2\pi y)$.
\end{example}

\begin{table}[htb!]
\begin{center}
\caption{Time accuracy for $\alpha=0.6,\,\sigma=0.4.$}\label{table:accuracy test alpha 06} \vspace*{0.3pt}
\def\temptablewidth{1.0\textwidth}
{\rule{\temptablewidth}{0.5pt}}
\begin{tabular*}{\temptablewidth}{@{\extracolsep{\fill}}cccccccccc}
\multirow{2}{*}{$N$} &\multirow{2}{*}{$\tau$} &\multicolumn{2}{c}{$\gamma=2$} &\multirow{2}{*}{$\tau$} &\multicolumn{2}{c}{$\gamma=4=\gamma_{\text{opt}}$} &\multirow{2}{*}{$\tau$}&\multicolumn{2}{c}{$\gamma=5$} \\
             \cline{3-4}          \cline{6-7}         \cline{9-10}
         &          &$e(N)$   &Order &         &$e(N)$   &Order &          &$e(N)$    &Order\\
\midrule
  200     &1.39e-02	 &8.33e-03 &$-$   &1.76e-02 &7.05e-04 &$-$   &1.68e-02	&3.82e-04 &$-$\\
  400     &7.26e-03	 &4.81e-03 &0.84  &9.24e-03 &2.47e-04 &1.63  &8.62e-03	&1.13e-04 &1.83\\
  800     &3.66e-03	 &2.77e-03 &0.81  &4.33e-03 &8.44e-05 &1.41  &4.54e-03	&3.24e-05 &1.95\\
  1600    &1.94e-03	 &1.59e-03 &0.87  &2.15e-03 &2.87e-05 &1.55  &2.20e-03	&9.12e-06 &1.75\\
  3200    &9.19e-04	 &9.13e-04 &0.74  &1.10e-03 &9.61e-06 &1.63  &1.13e-03	&2.99e-06 &1.68\\
\midrule
\multicolumn{3}{l}{$\min\{1+\alpha, \gamma\sigma\}$}   &0.80 & & &1.60 & & &1.60\\
\end{tabular*}
{\rule{\temptablewidth}{0.5pt}}
\end{center}
\end{table}	

\begin{table}[htb!]
\begin{center}
\caption{Time accuracy for $\alpha=0.8,\,\sigma=0.6.$}\label{table:accuracy test alpha 08} \vspace*{0.3pt}
\def\temptablewidth{1.0\textwidth}
{\rule{\temptablewidth}{0.5pt}}
\begin{tabular*}{\temptablewidth}{@{\extracolsep{\fill}}cccccccccc}
\multirow{2}{*}{$N$} &\multirow{2}{*}{$\tau$} &\multicolumn{2}{c}{$\gamma=2$} &\multirow{2}{*}{$\tau$} &\multicolumn{2}{c}{$\gamma=3=\gamma_{\text{opt}}$} &\multirow{2}{*}{$\tau$}&\multicolumn{2}{c}{$\gamma=4$} \\
             \cline{3-4}          \cline{6-7}         \cline{9-10}
         &          &$e(N)$   &Order &         &$e(N)$   &Order &          &$e(N)$    &Order\\
\midrule
  200     &1.42e-02	&6.97e-04 &$-$   &1.58e-02 &1.31e-04 &$-$   &1.74e-02  &6.44e-05  &$-$\\
  400     &8.02e-03	&3.07e-04 &1.44  &8.73e-03 &4.07e-05 &1.97  &8.62e-03  &1.71e-05  &1.88\\
  800    &3.73e-03	&1.34e-04 &1.08  &4.13e-03 &1.25e-05 &1.58  &4.46e-03  &4.39e-06  &2.07\\
  1600    &1.93e-03	&5.86e-05 &1.26  &2.07e-03 &3.77e-06 &1.73  &2.10e-03  &1.13e-06  &1.80\\
  3200    &9.53e-04	&2.55e-05 &1.18  &1.06e-03 &1.13e-06 &1.82  &1.13e-03  &4.15e-07  &1.63\\
\midrule
\multicolumn{3}{l}{$\min\{1+\alpha, \gamma\sigma\}$}   &1.20 & & &1.80 & & &1.80\\
\end{tabular*}
{\rule{\temptablewidth}{0.5pt}}
\end{center}
\end{table}	

To resolve the initial singularity, we split the time interval $[0,T]$ into two parts $[0, T_{0}]$ and $[T_{0}, T]$
with total $N$ subintervals.
A graded mesh $t_{k}=T_{0}(k/N_0)^{\gamma}$ is employed with
$T_0=\min\{1/\gamma,T\}$ and $N_0=\lceil \frac{N}{T+1-\gamma^{-1}}\rceil$
in the first part $[0, T_{0}]$.
In the remainder interval $[T_{0},T]$, we use random step sizes
$$\tau_{N_{0}+k}:=\frac{(T-T_{0})\epsilon_{k}}{\sum_{k=1}^{N-N_0}\epsilon_{k}}\quad \text{for $1\leq k\leq N-N_0$,}$$
where $\epsilon_{k}$ are uniformly distributed random numbers inside $(0,1)$.

We focus on the time accuracy of the modified Crank-Nicolson scheme
\eqref{eq: L1rCN scheme1}-\eqref{eq: L1rCN scheme2}.
Always, the spatial domain $\Omega=(0,1)^2$ is uniformly discretized by using
$512 \times 512$ grids such that the temporal error dominates.
We record the maximum norm error
$e(N):=\max_{1\leq{n}\leq{N}}\|u(t_n)-u^n\|_{\infty}$
in each run
and evaluate the convergence order by
$$\text{Order}\approx\frac{\log\bra{e(N)/e(2N)}}{\log\bra{\tau(N)/\tau(2N)}},$$
where $\tau(N)$ denotes the maximum time-step size for total $N$ subintervals.
We run the new scheme by considering the following two cases:
\begin{enumerate}
\item[(a)] The fractional order $\alpha=0.6$ and regularity parameter $\sigma=0.4$ with
mesh parameters  $\gamma=2,4,5$, respectively (see Table \ref{table:accuracy test alpha 06});
\item[(b)] The fractional order $\alpha=0.8$ and regularity parameter $\sigma=0.6$ with
mesh parameters $\gamma=2,3,4$, respectively (see Table \ref{table:accuracy test alpha 08}).
\end{enumerate}

From  Tables \ref{table:accuracy test alpha 06}
and \ref{table:accuracy test alpha 08},
one can observe that an optimal rate
$\mathcal{O}\bra{\tau^{\{\gamma\sigma,1+\alpha\}}}$
is achieved when the grading parameter
$\gamma\ge\gamma_{\text{opt}}=\max\{1,(1+\alpha)/\sigma\}$.
As noticed, the error analysis in \cite{Mustapha:2011,MustaphaAlMutawa:2012} is only suited for the graded meshes.
Thus there is still a gap between the numerical evidences and the theoretical
verification of convergence rates on a general class of nonuniform time meshes.

\subsection{Discrete maximum bound principle}

\begin{figure}[htb!]
\centering
\includegraphics[width=2.0in]{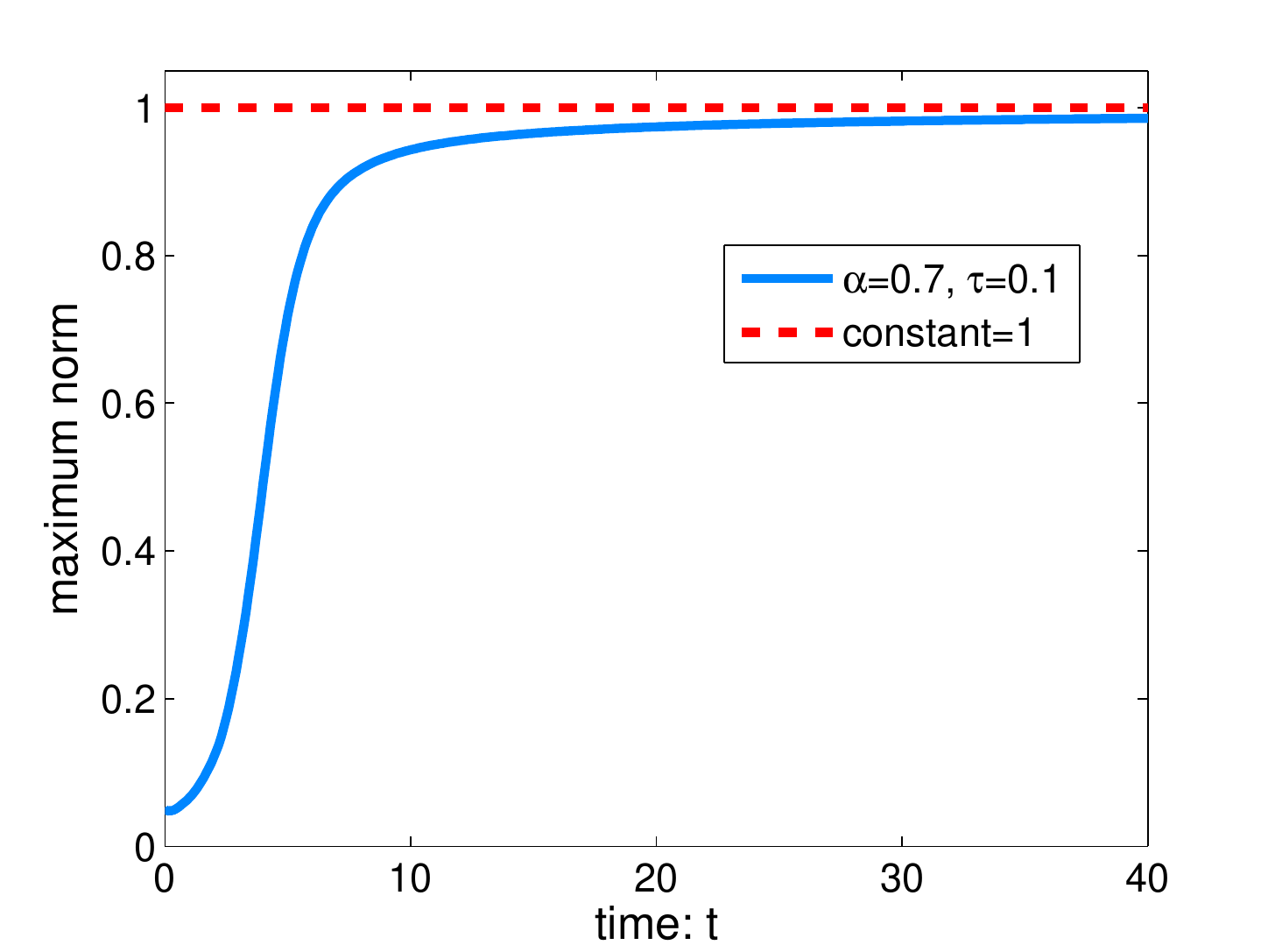}
\includegraphics[width=2.0in]{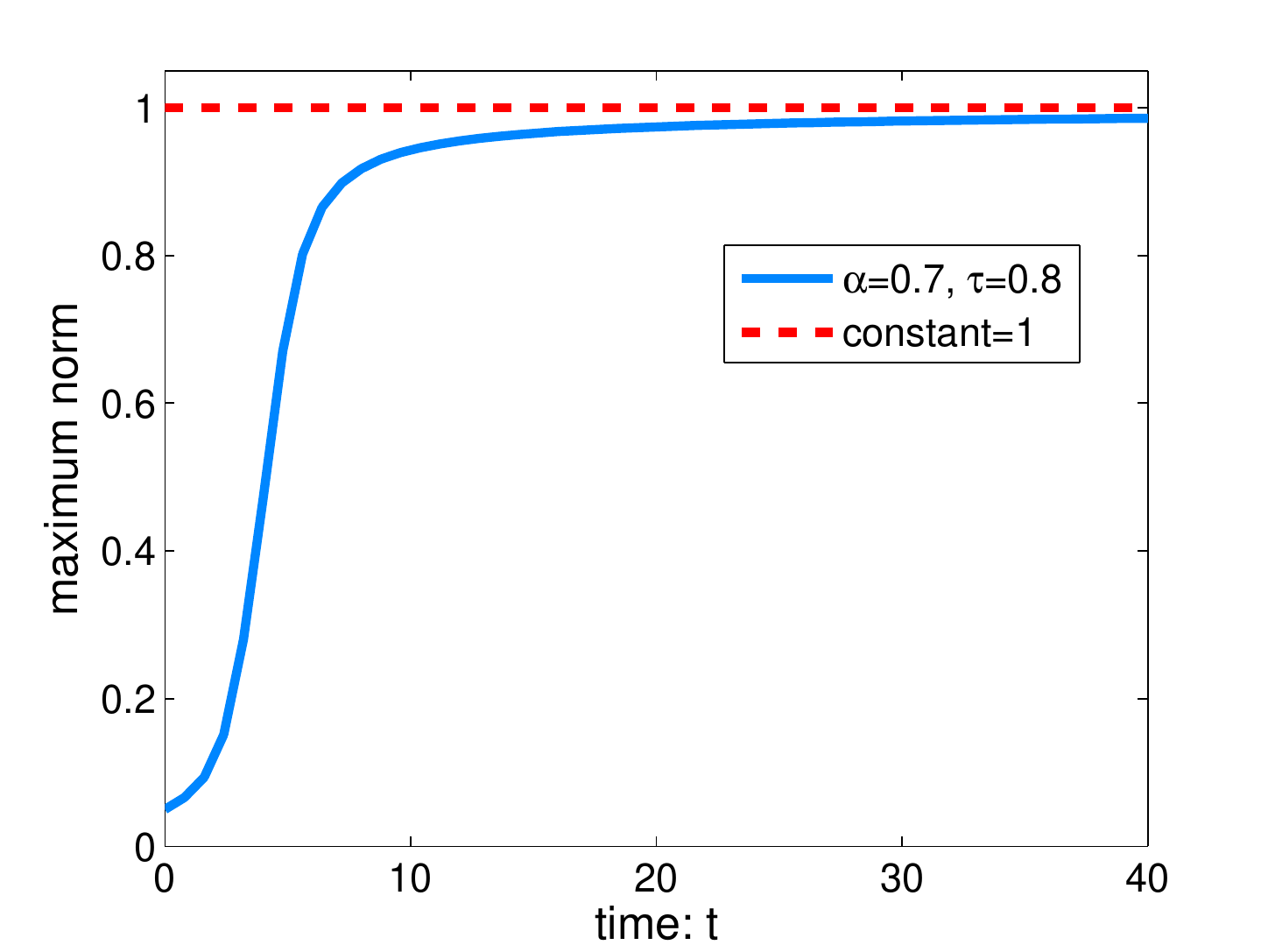}
\includegraphics[width=2.0in]{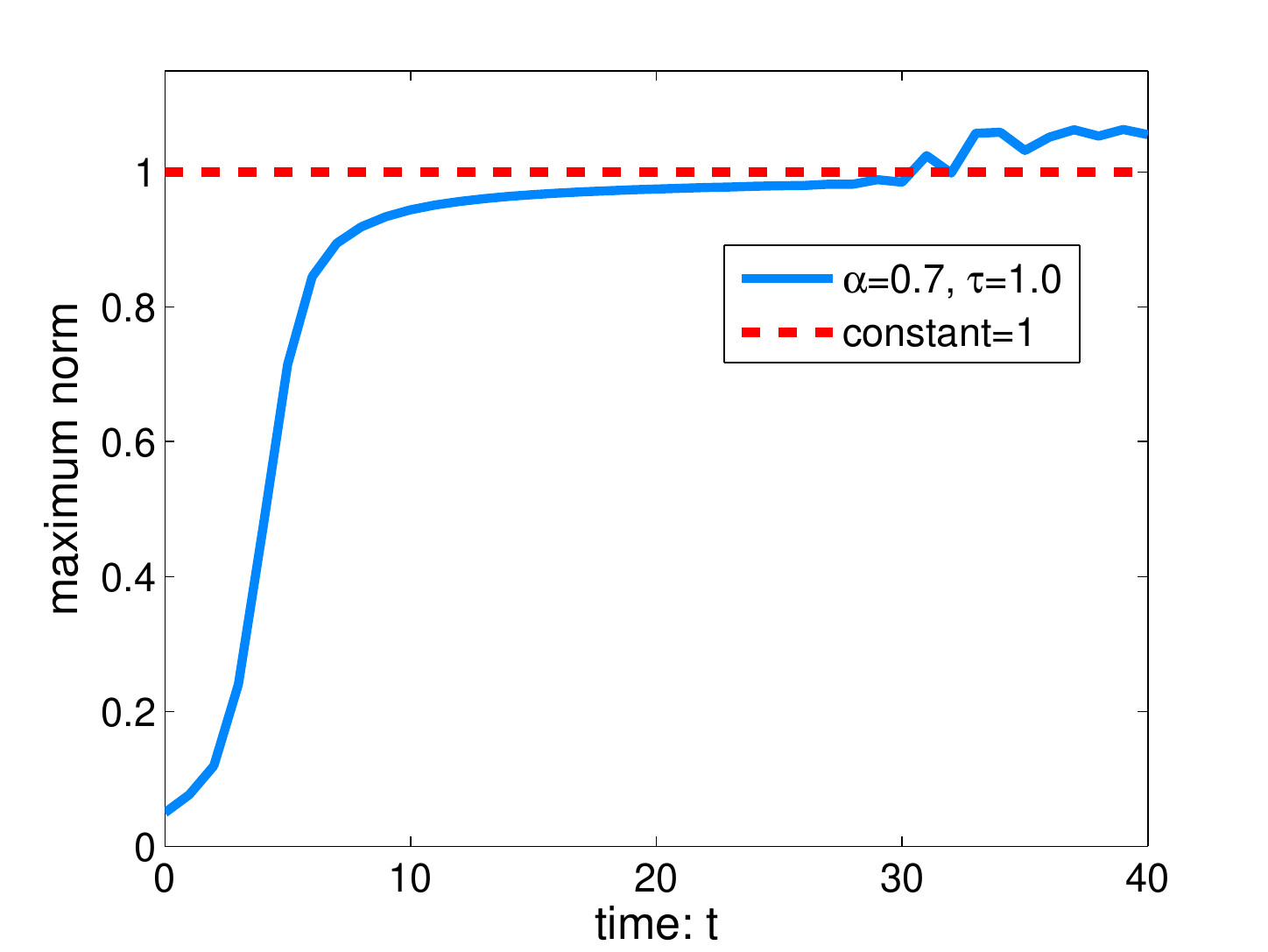}\\
\includegraphics[width=2.0in]{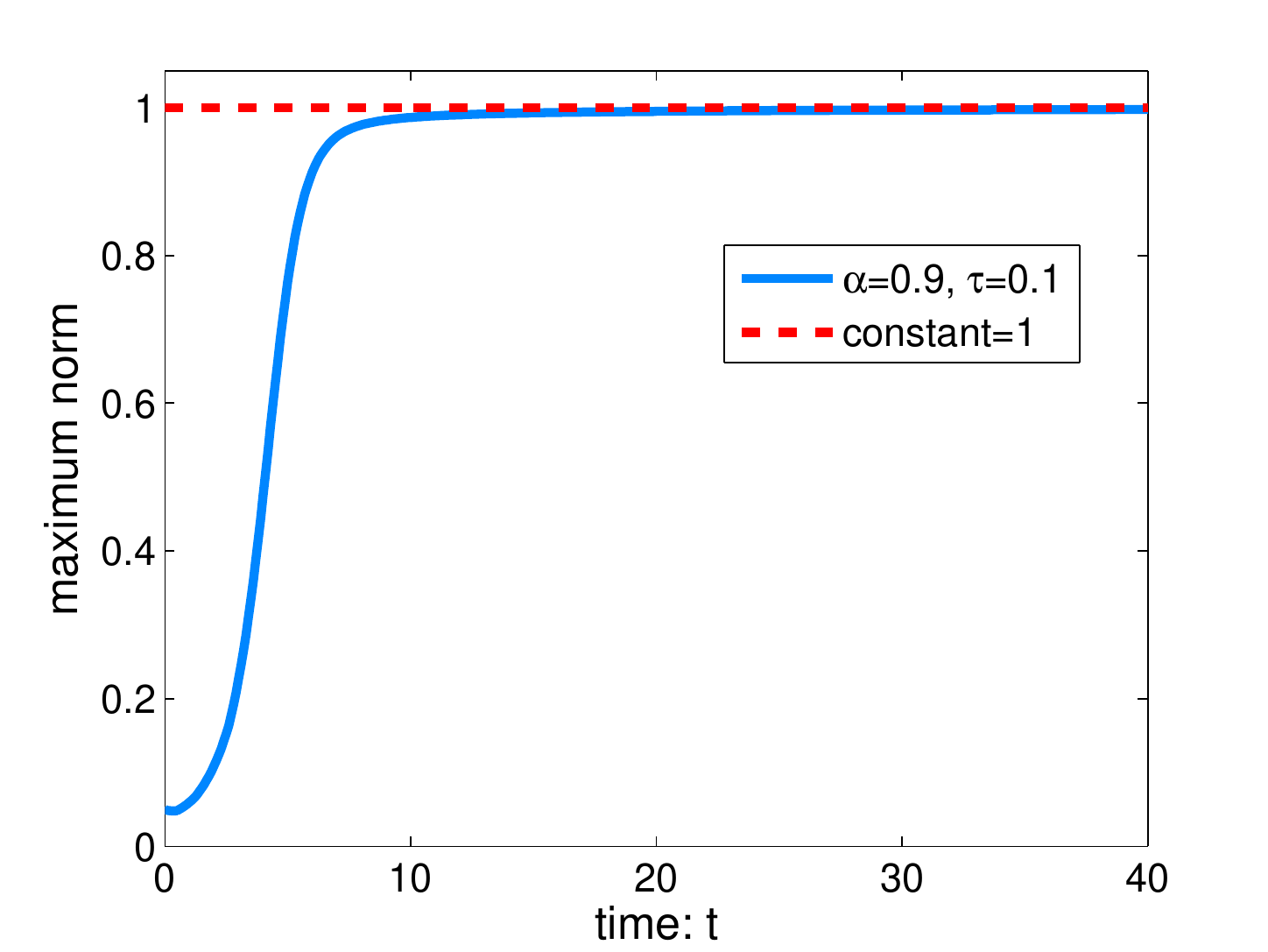}
\includegraphics[width=2.0in]{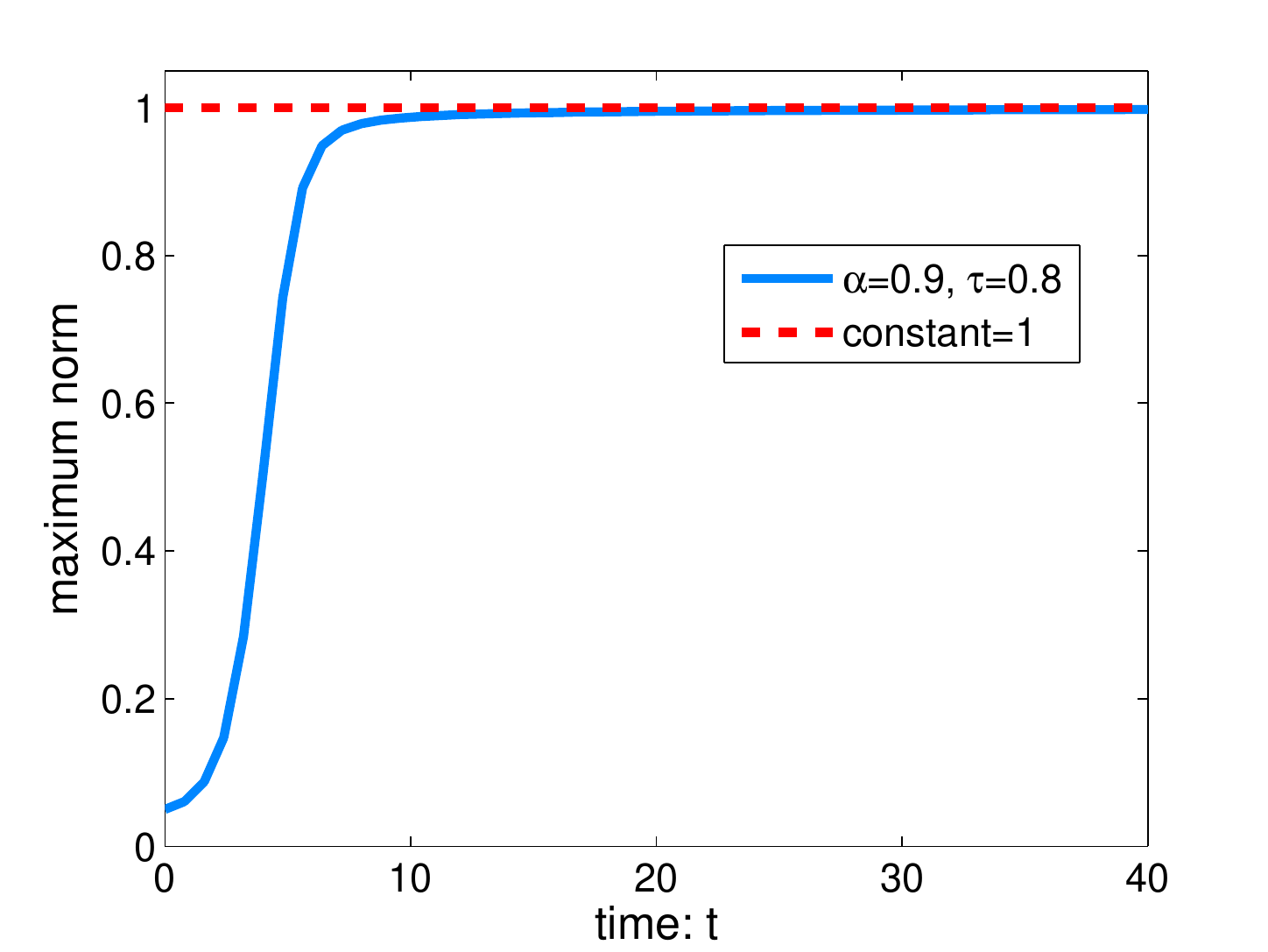}
\includegraphics[width=2.0in]{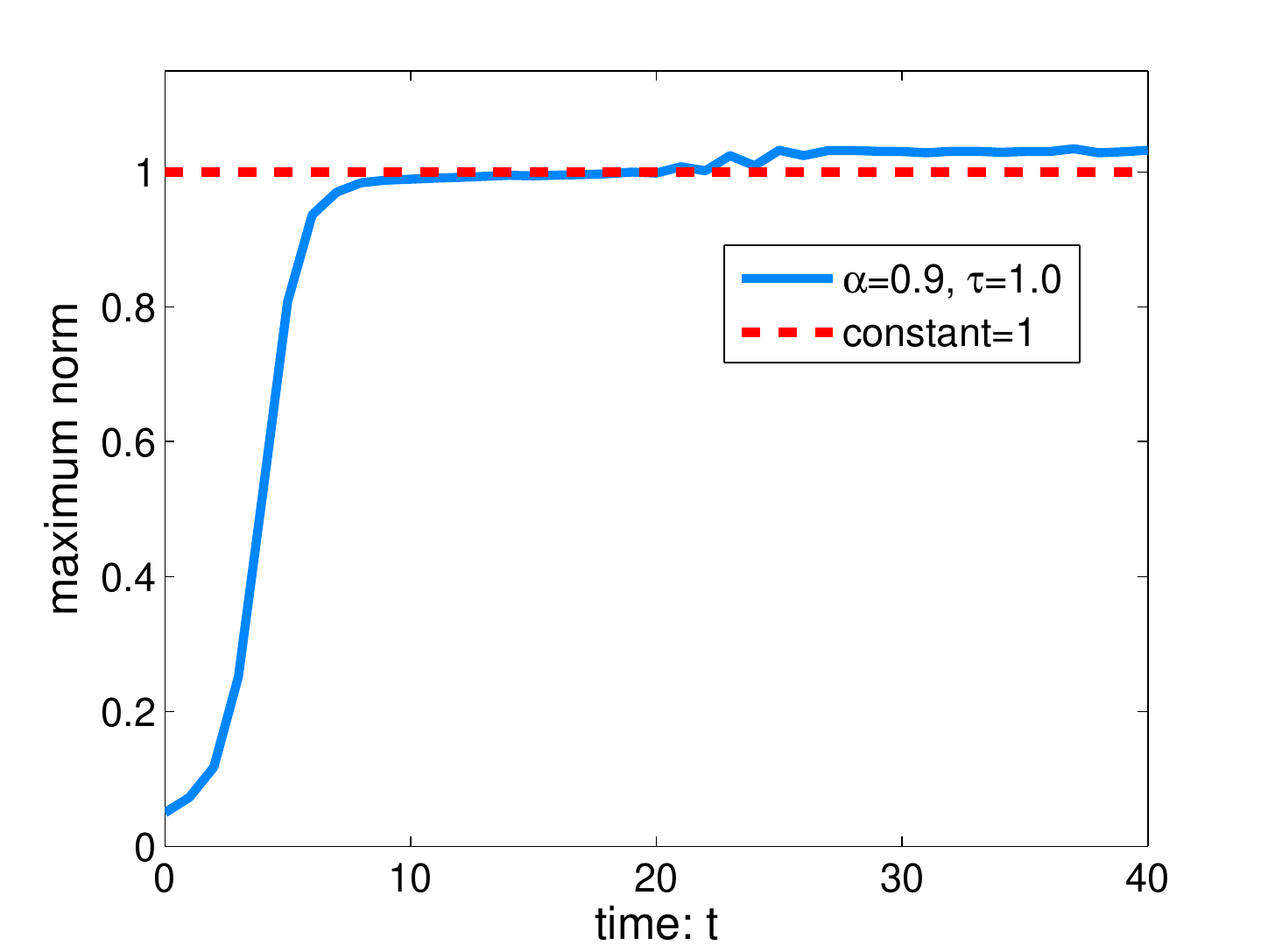}\\
\caption{The maximum norm of discrete solutions  for the fractional orders $\alpha=0.7,0.9$
  (from top to bottom) with three different time-step sizes $\tau=0.1,0.8,1.0$ (from left to right).}
\label{figs:maximum principle}
\end{figure}

We now verify the discrete  maximum bound principle.
For the fractional orders $\alpha=0.7$ and $0.9$,
we run the numerical scheme \eqref{eq: L1rCN scheme1}-\eqref{eq: L1rCN scheme2}
with the  random initial data $u_0(\mathbf{x})=\mathrm{rand}(\mathbf{x})$ until $T=40$ on different uniform meshes.
Figure \ref{figs:maximum principle} plots the maximum norm for two fractional order $\alpha=0.7,0.9$
with three different time-step size $\tau=0.1, 0.8, 1.0$.
These results suggest that the time-step restriction
\eqref{Constraint: Dis_Max Time-Steps} is only sufficient to
ensure the maximum maximum principle.
Actually, the step-size restriction \eqref{Constraint: Dis_Max Time-Steps}
requires the maximal step size $\tau\le0.14$ for the fractional order $\alpha=0.7$
and requires $\tau\le0.36$ as the fractional order $\alpha=0.9$.

\subsection{Initial singularity and graded meshes}

\begin{example}\label{exam:coarsen dynamic}
Consider the time-fractional
Allen-Cahn equation \eqref{Cont: TFAC} on the physical domain $(0,2\pi)^2$
with the model parameter $\varepsilon=0.05$.
The initial condition is taken as  $u_0(\mathbf{x})=\mathrm{rand}(\mathbf{x})$,
where $\mathrm{rand}(\mathbf{x})$ is uniformly distributed random number
varying from $-0.001$ to $0.001$ to each grid points. Always, a $128 \times 128$ uniform spatial mesh
is used to cover the domain $(0,2\pi)^2$.
\end{example}

\begin{figure}[htb!]
\centering
\includegraphics[width=3.0in,height=2.0in]{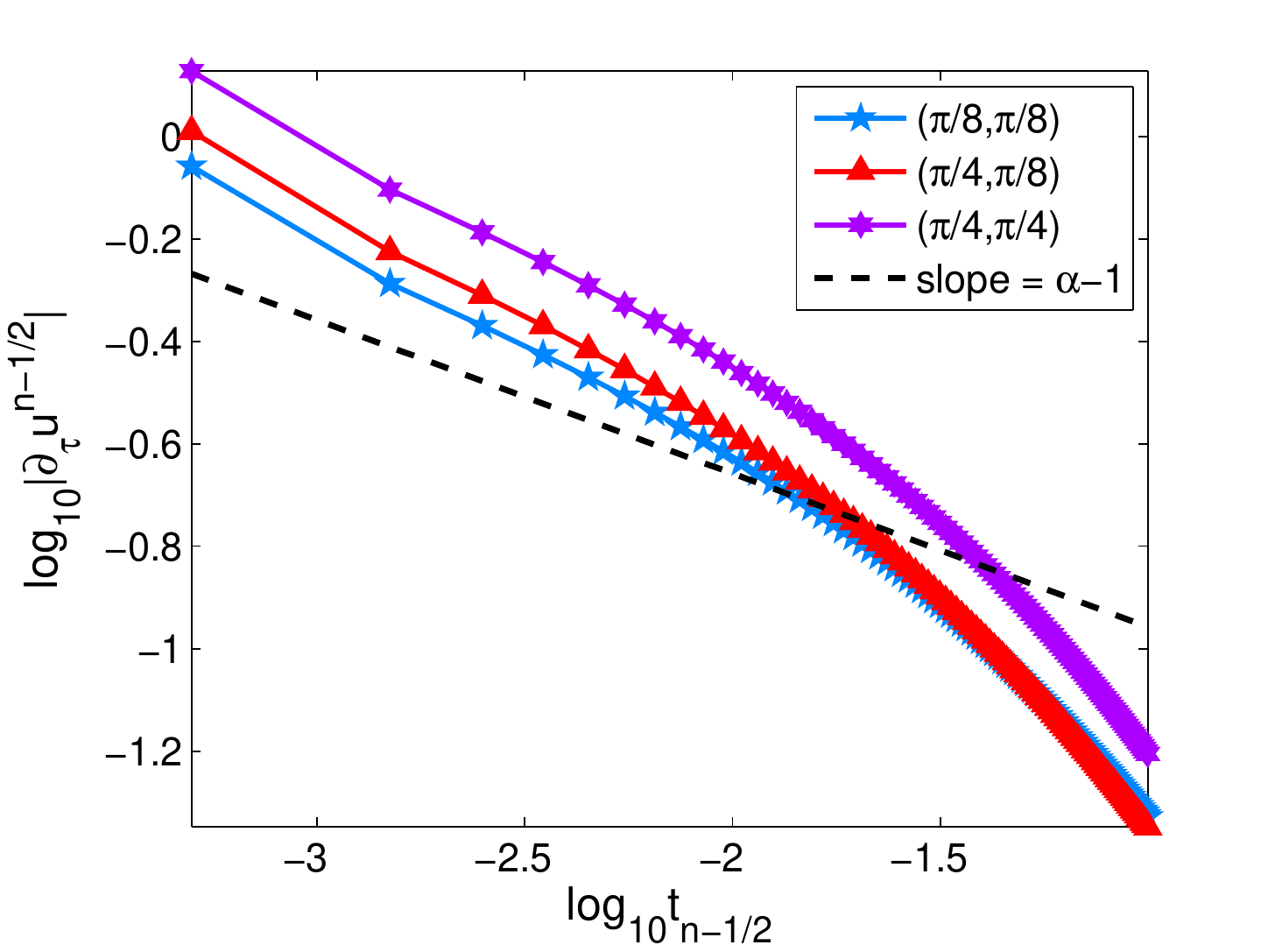}
\includegraphics[width=3.0in,height=2.0in]{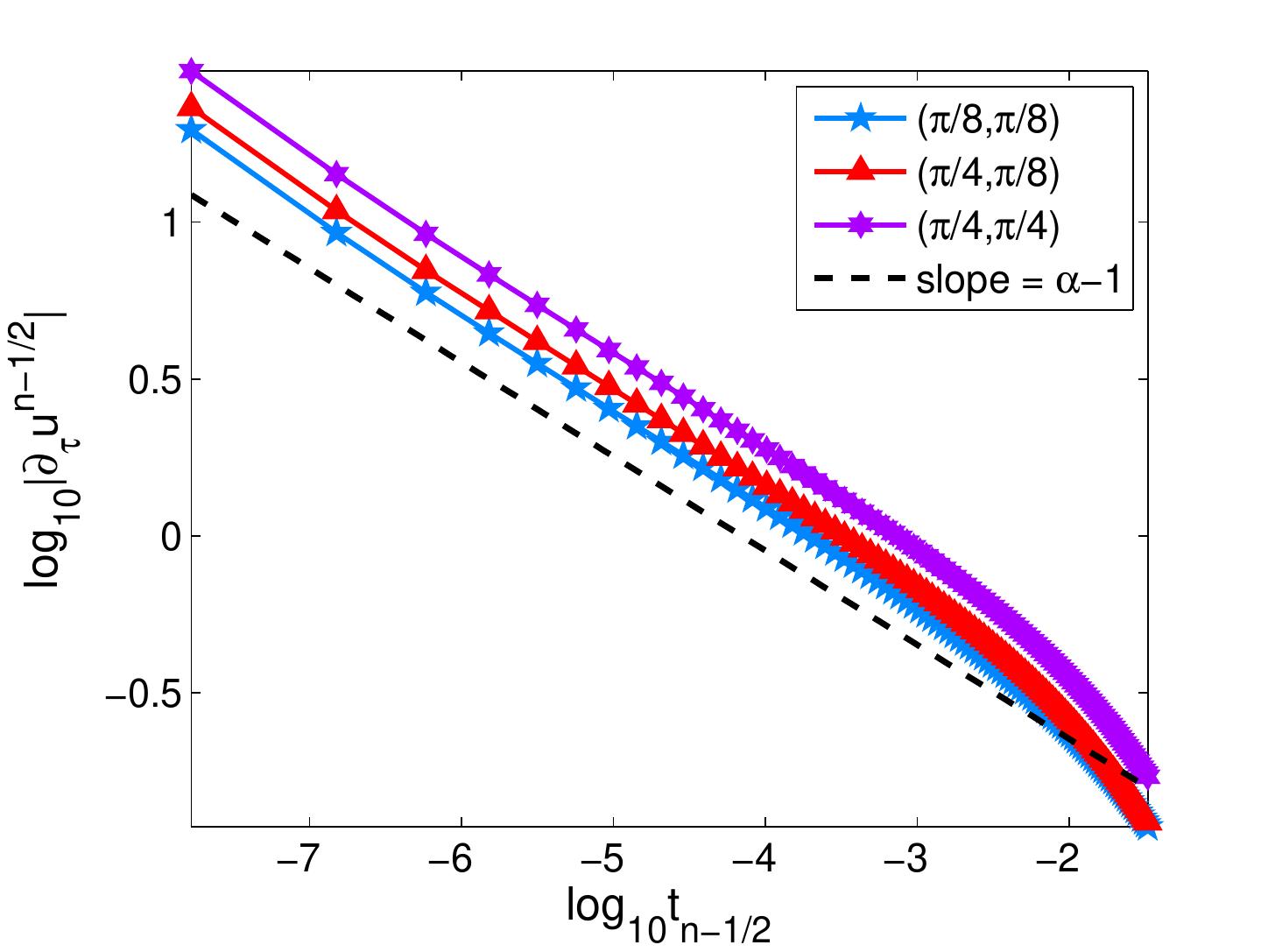}
\caption{The log-log plot of $\partial_{\tau}u^{n-\frac12}$ versus time for \eqref{Cont: TFAC} with $\alpha=0.7$
and grading parameters $\gamma=1,\,3$ (from left to right). The legends refer to the spatial positions.}
\label{figs:singu RL AC alpha 07}
\end{figure}

We run the scheme \eqref{eq: L1rCN scheme1}-\eqref{eq: L1rCN scheme2} with fractional order $\alpha=0.7$,
$T=1/\gamma$ and $u_0(\mathbf{x})=\mathrm{rand}(\mathbf{x})$.
Figure \ref{figs:singu RL AC alpha 07}
depicts the discrete time derivative $\partial_{\tau}u^{n-\frac12}$
near $t=0$ on the graded mesh $t_n=T(n/N)^{\gamma}$ for two grading parameters $\gamma=1,3$.
It is seen that
\[
\log|u_{t}(\mathbf{x},t)|
\approx(\alpha-1)\log(t)+C(\mathbf{x})\;\;\text{such that $u_t=\mathcal{O}(t^{\alpha-1})$}\quad\text{as $t\rightarrow0$,}
\]
and the initial singularity can be resolved
by concentrating more grids near initial time.

\subsection{Adaptive time stepping}

\begin{figure}[htb!]
\centering
\includegraphics[width=2.0in]{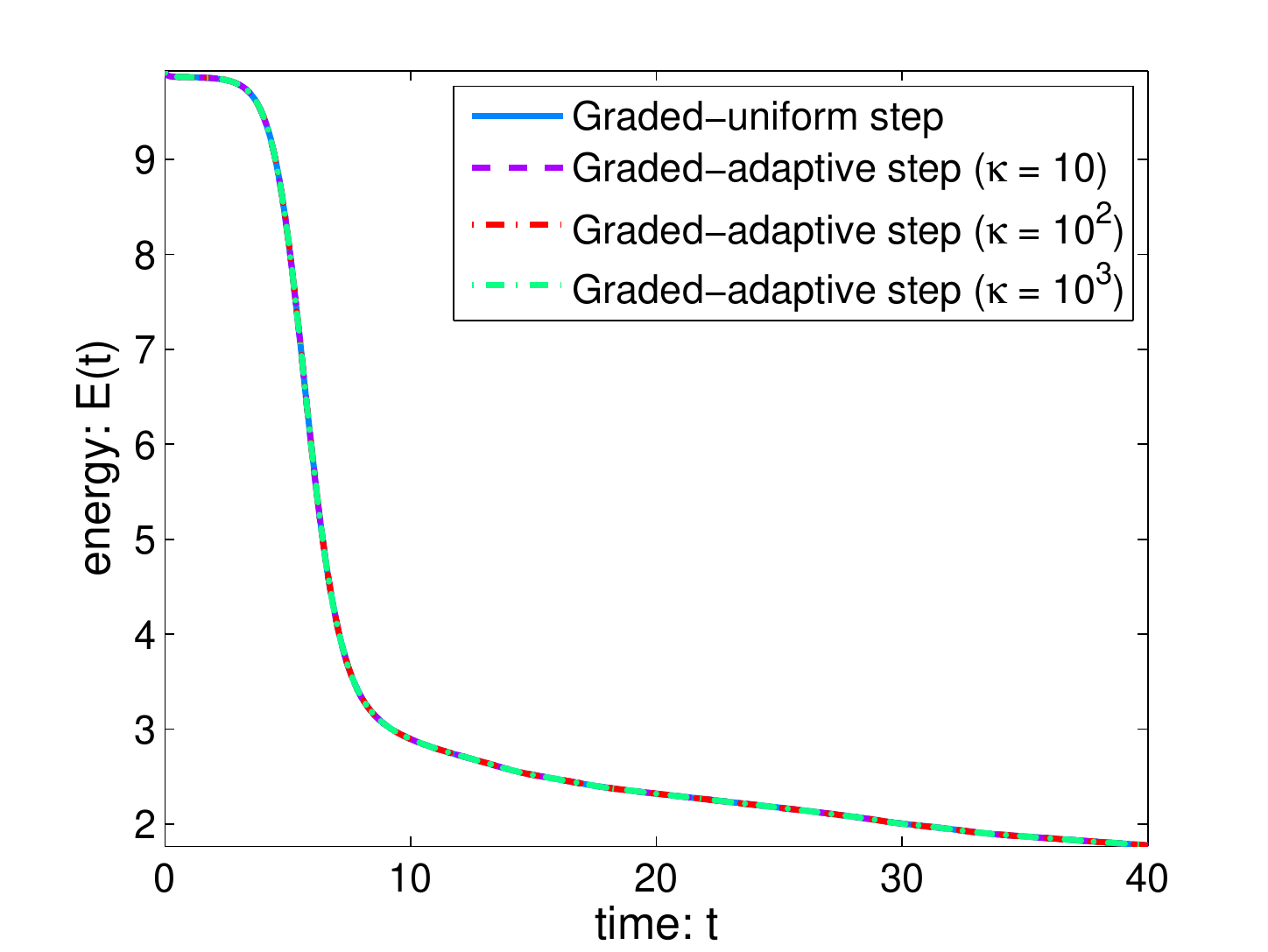}
\includegraphics[width=2.0in]{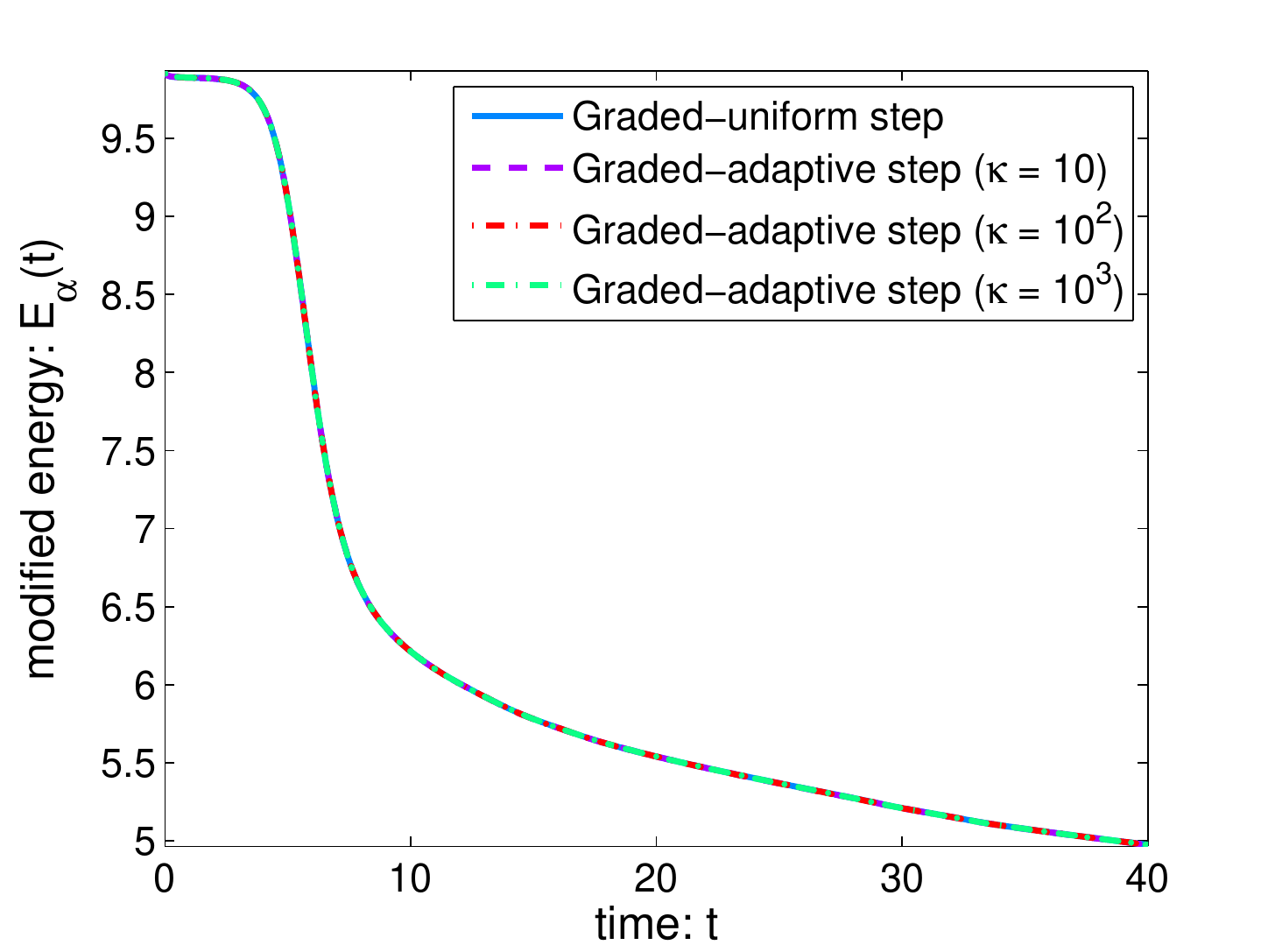}
\includegraphics[width=2.0in]{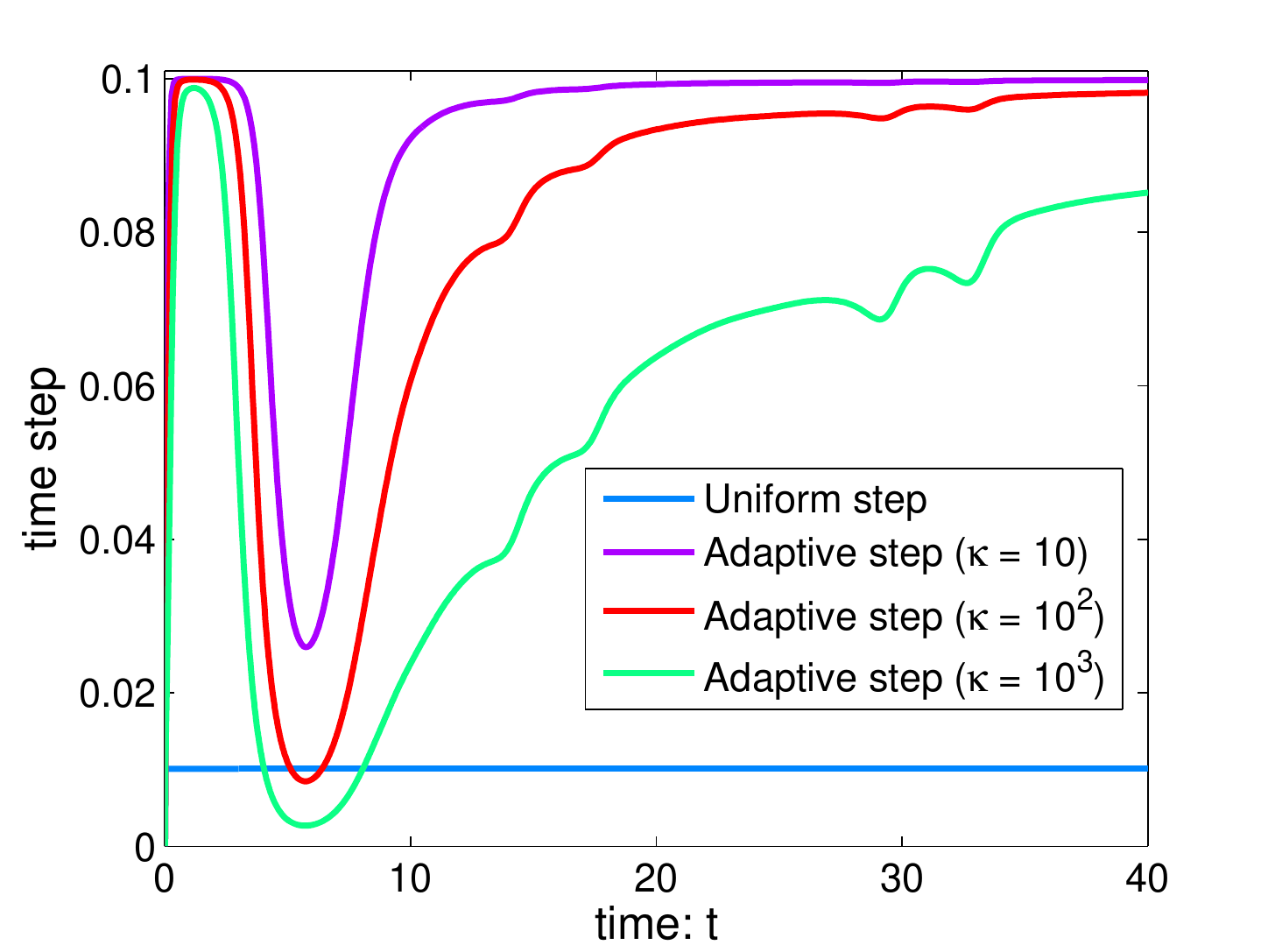}
\caption{The energies $E(t)$, $E_\alpha(t)$ and adaptive steps of Example \ref{exam:coarsen dynamic}.}
\label{figs:compar time steps}
\end{figure}

\begin{table}[htb!]
\begin{center}
\caption{Comparisons of CPU time (in seconds) and time steps.}\label{table:compar CPU time}
\vspace*{0.3pt}
\def\temptablewidth{0.9\textwidth}
{\rule{\temptablewidth}{0.5pt}}
\begin{tabular*}{\temptablewidth}{@{\extracolsep{\fill}}c|ccccc}
  Adaptive parameter & $\kappa=10$ & $\kappa=10^2$ & $\kappa=10^3$ &uniform mesh\\
  \midrule
  CPU time  &41.167               &61.596                 &136.787                &321.830\\
  Time steps     &507                  &769                    &1734                   &4000\\
\end{tabular*}
{\rule{\temptablewidth}{0.5pt}}
\end{center}
\end{table}	
In order to resolve the dynamic evolutions involving multiple time scales
and reduce the computation cost in long time simulations,
we next present an adaptive time-stepping strategy \cite{QiaoZhangTang:2011}
with the following adaptive step criterion based on the energy variation,
\begin{align*}
\tau_{ada}
=\max\Bigg\{\tau_{\min},
\frac{\tau_{\max}}{\sqrt{1+\kappa\abs{E_\alpha^{\prime}(t)}^{2}}}\Bigg\},
\end{align*}
where $E_\alpha$ is the modified energy \eqref{Cont: pseudo-local energy},
$\tau_{\max},\tau_{\min}$ are
the predetermined maximum and minimum time steps, respectively.
The parameter $\kappa$ is chosen to adjust the level of adaptivity.
In our computations,  the time interval $[0,T]$ is divided into $[0,T_0]$ and $[T_0,T]$.
We choose the graded mesh $t_k=T_0(k/N_0)^\gamma$ with $T_0=0.01$, $N_0=30$ and $\gamma=3$ in the starting cell $[0,T_0]$.
The remainder $[T_0,T]$ is tested by two types of time meshes:
\begin{description}
\item[(Graded-uniform mesh)] Uniform step size $\tau=0.01$;
  \item[(Graded-adaptive mesh)] Adaptive time-stepping with $\tau_{\max}=10^{-1}$ and $\tau_{\min}=10^{-3}$.
\end{description}
\begin{figure}[htb!]
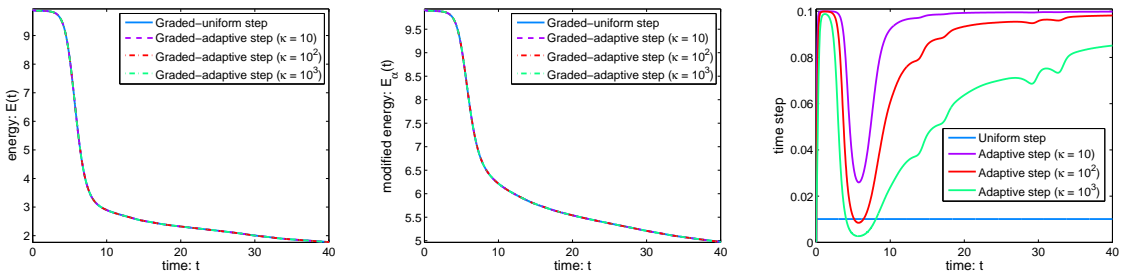

\centering
\includegraphics[width=2.0in]{compar_adap_energy_loc-eps-converted-to}
\includegraphics[width=2.0in]{compar_adap_energy_alpha-eps-converted-to}
\includegraphics[width=2.0in]{compar_adap_time_step-eps-converted-to}
\caption{The energies $E(t)$, $E_\alpha(t)$ and adaptive steps of Example \ref{exam:coarsen dynamic}.}
\label{figs:compar time steps}
\end{figure}

\begin{table}[htb!]
\begin{center}
\caption{Comparisons of CPU time (in seconds) and time steps.}\label{table:compar CPU time}
\vspace*{0.3pt}
\def\temptablewidth{0.9\textwidth}
{\rule{\temptablewidth}{0.5pt}}
\begin{tabular*}{\temptablewidth}{@{\extracolsep{\fill}}c|ccccc}
  Adaptive parameter & $\kappa=10$ & $\kappa=10^2$ & $\kappa=10^3$ &uniform mesh\\
  \midrule
  CPU time  &41.167               &61.596                 &136.787                &321.830\\
  Time steps     &507                  &769                    &1734                   &4000\\
\end{tabular*}
{\rule{\temptablewidth}{0.5pt}}
\end{center}
\end{table}	

\begin{figure}[htb!]
\centering
\includegraphics[width=1.47in]{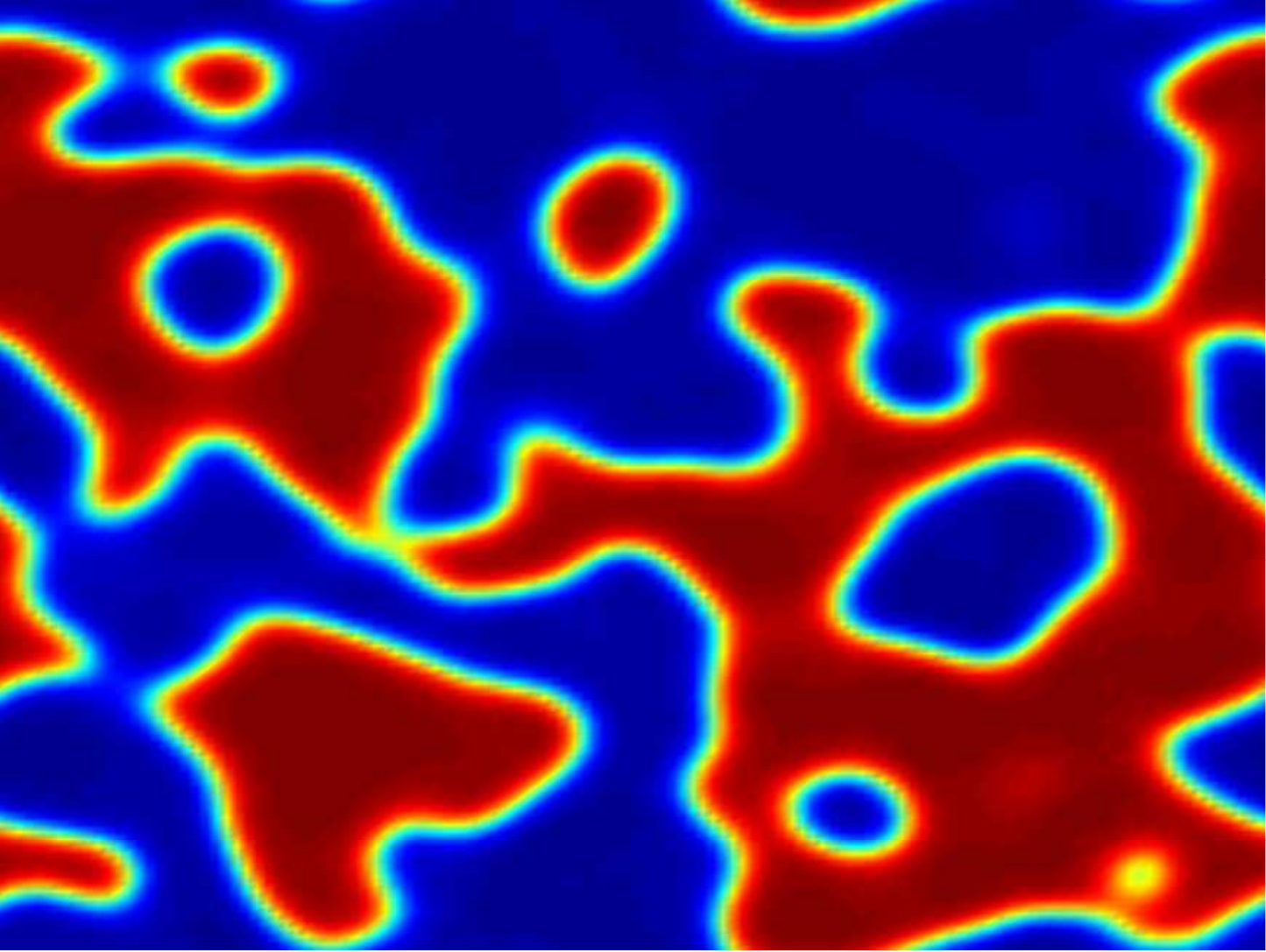}
\includegraphics[width=1.47in]{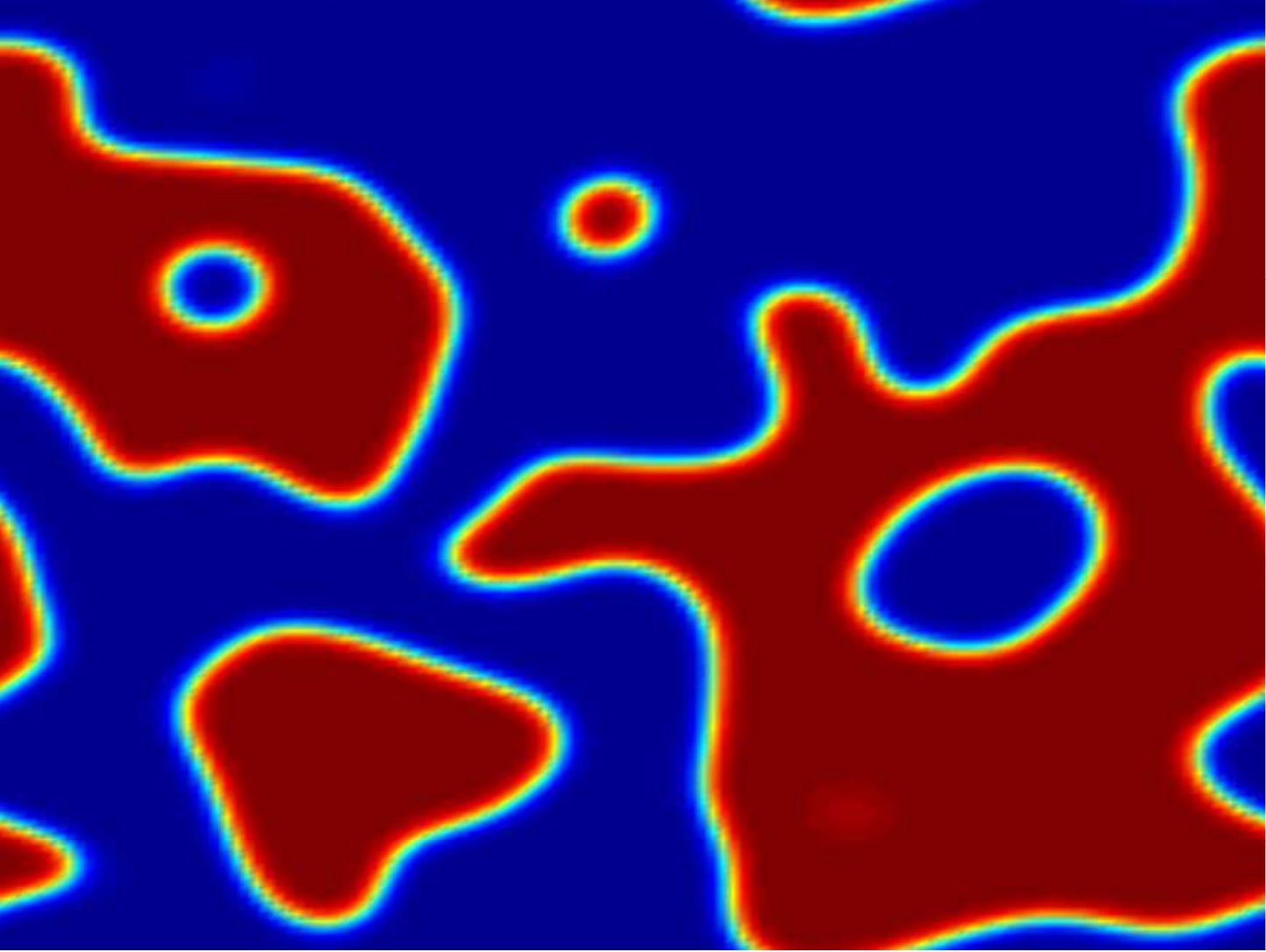}
\includegraphics[width=1.47in]{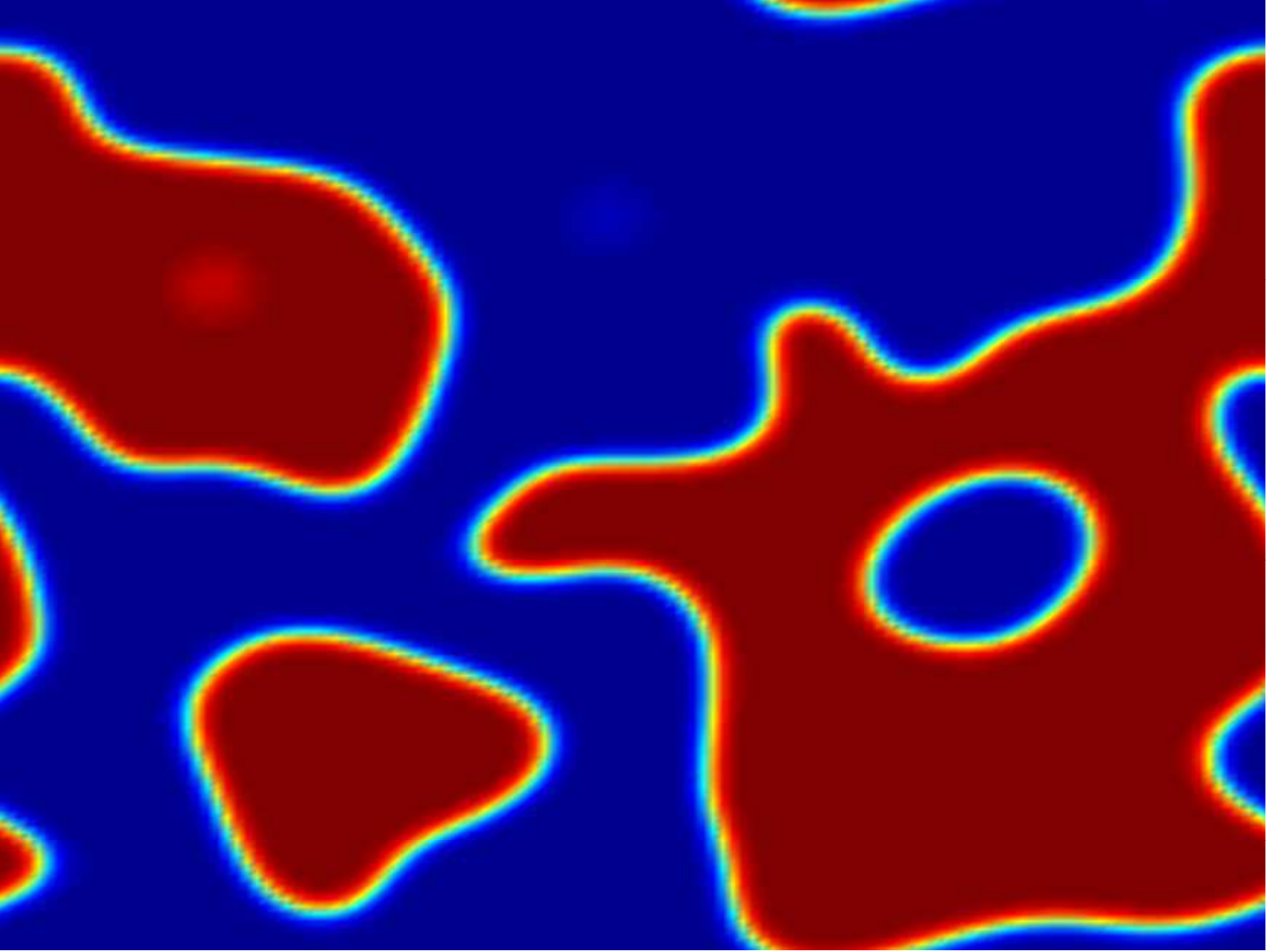}
\includegraphics[width=1.47in]{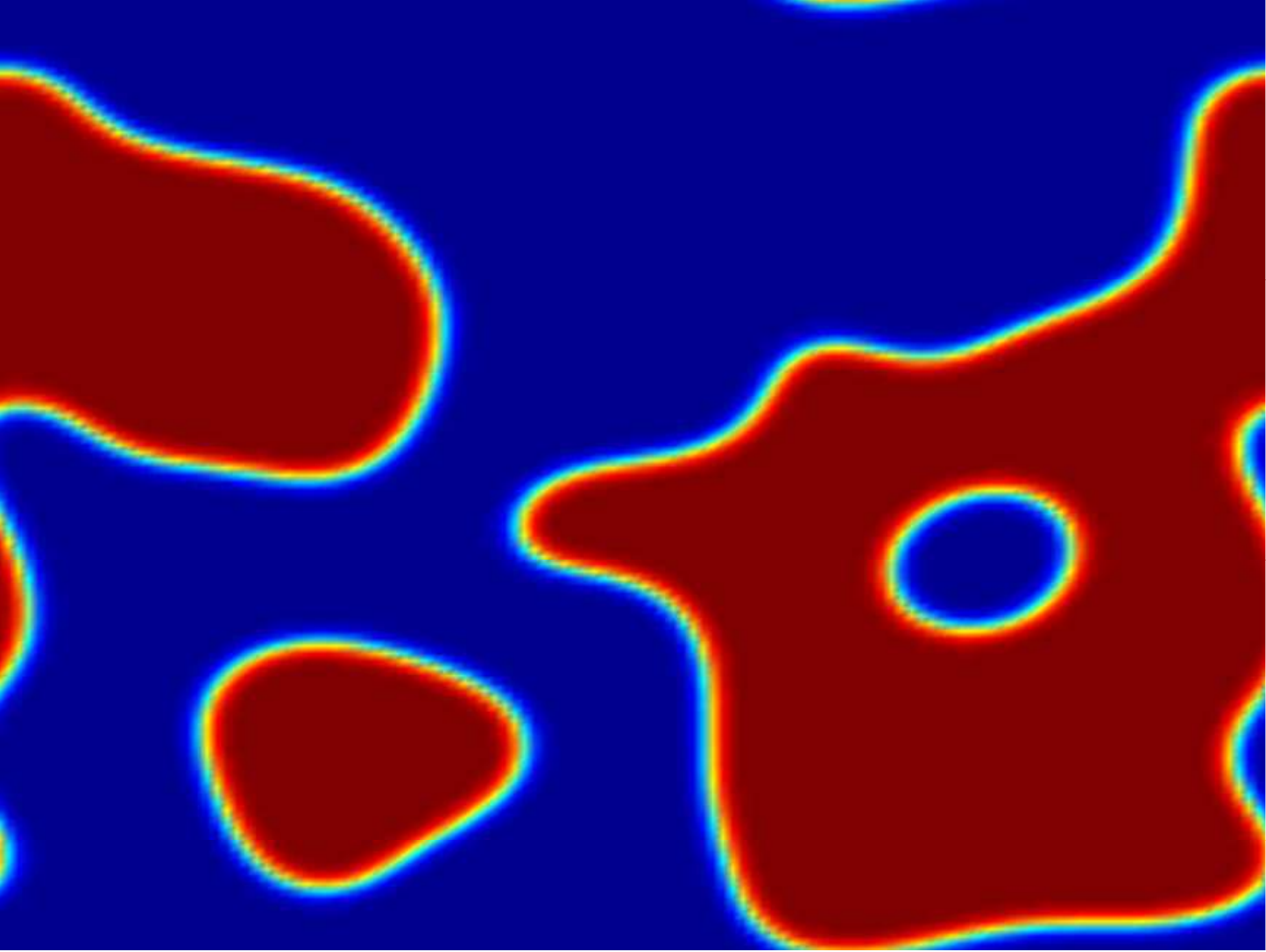}\\
\includegraphics[width=1.47in]{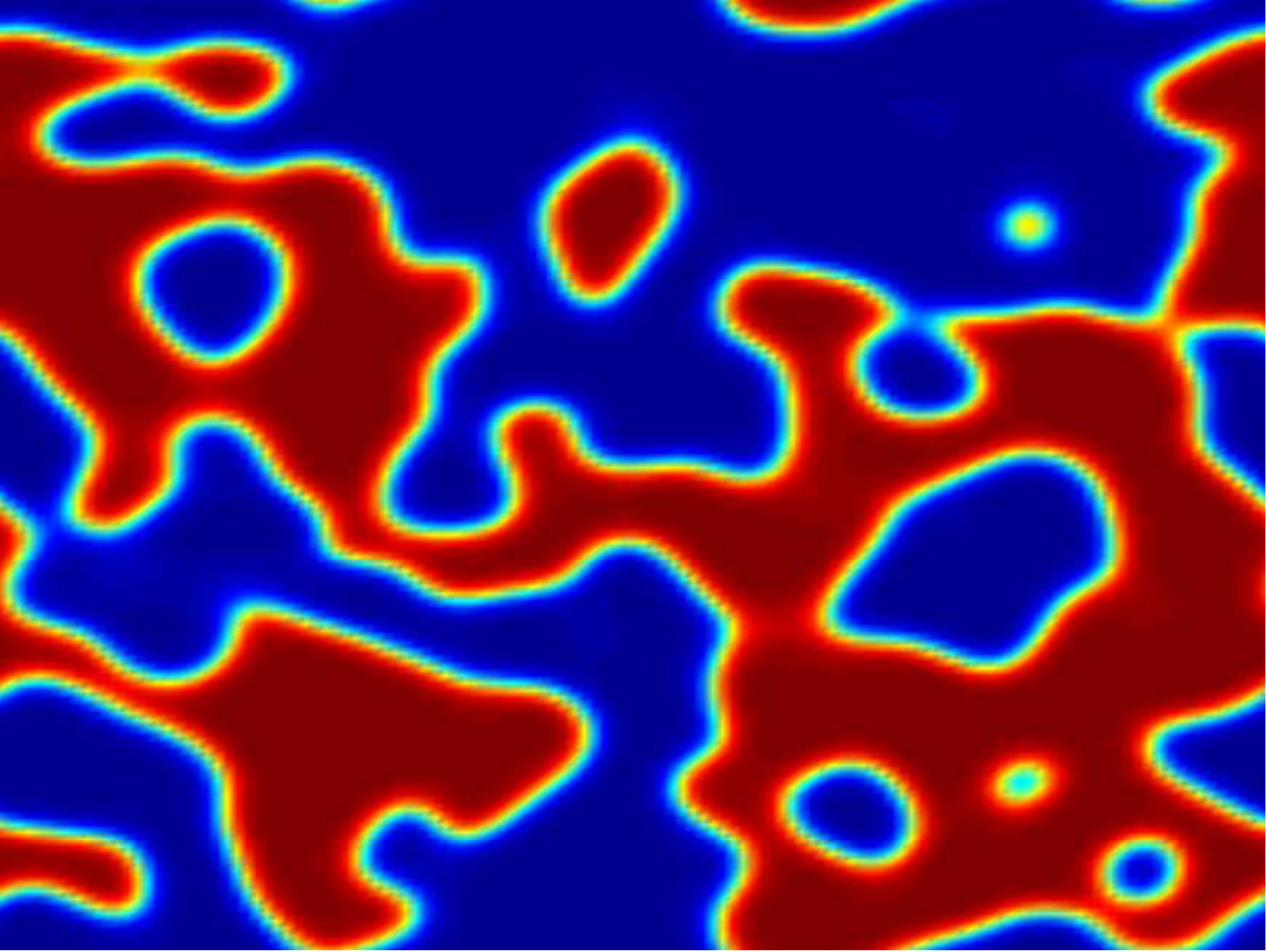}
\includegraphics[width=1.47in]{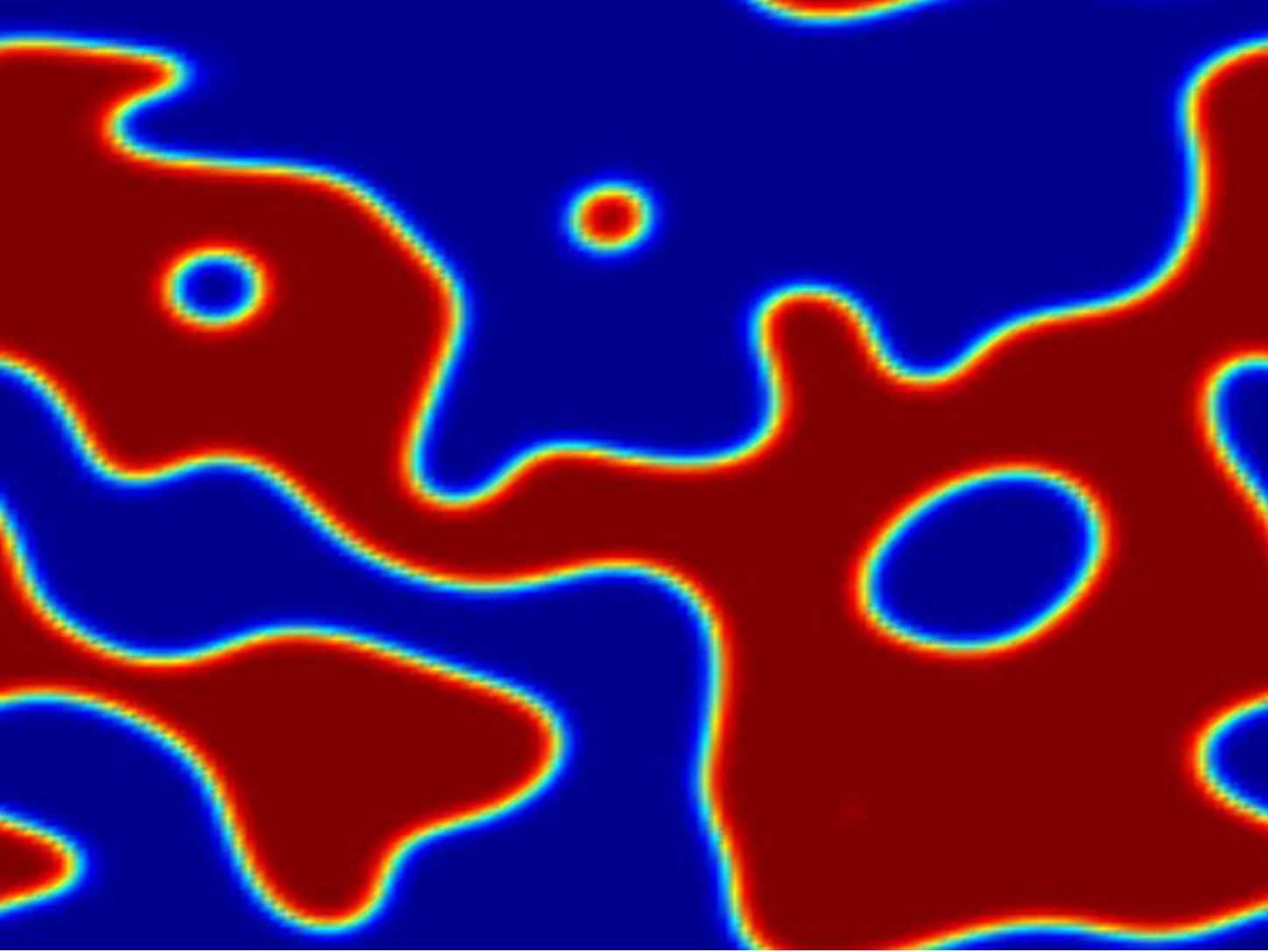}
\includegraphics[width=1.47in]{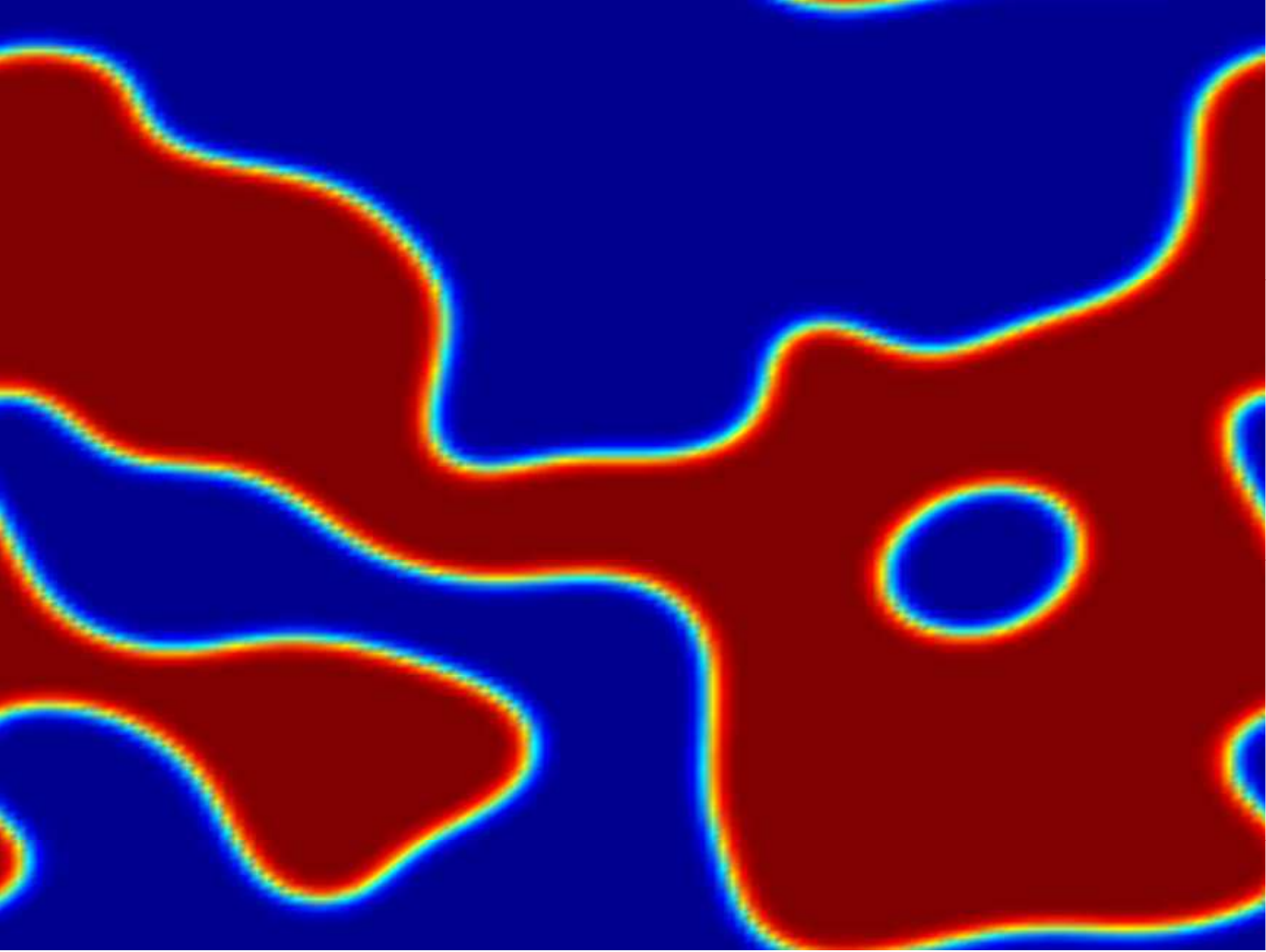}
\includegraphics[width=1.47in]{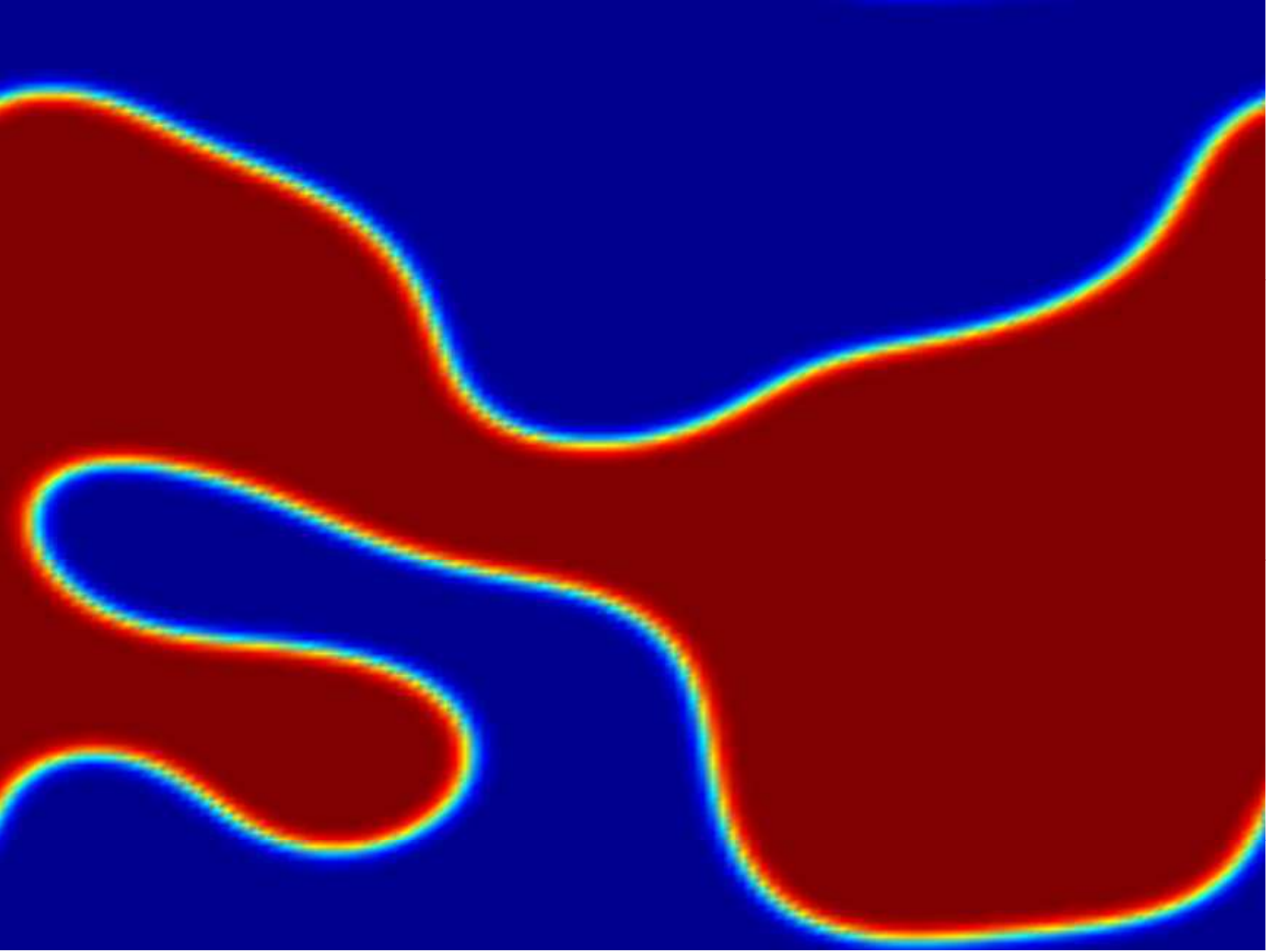}\\
\includegraphics[width=1.47in]{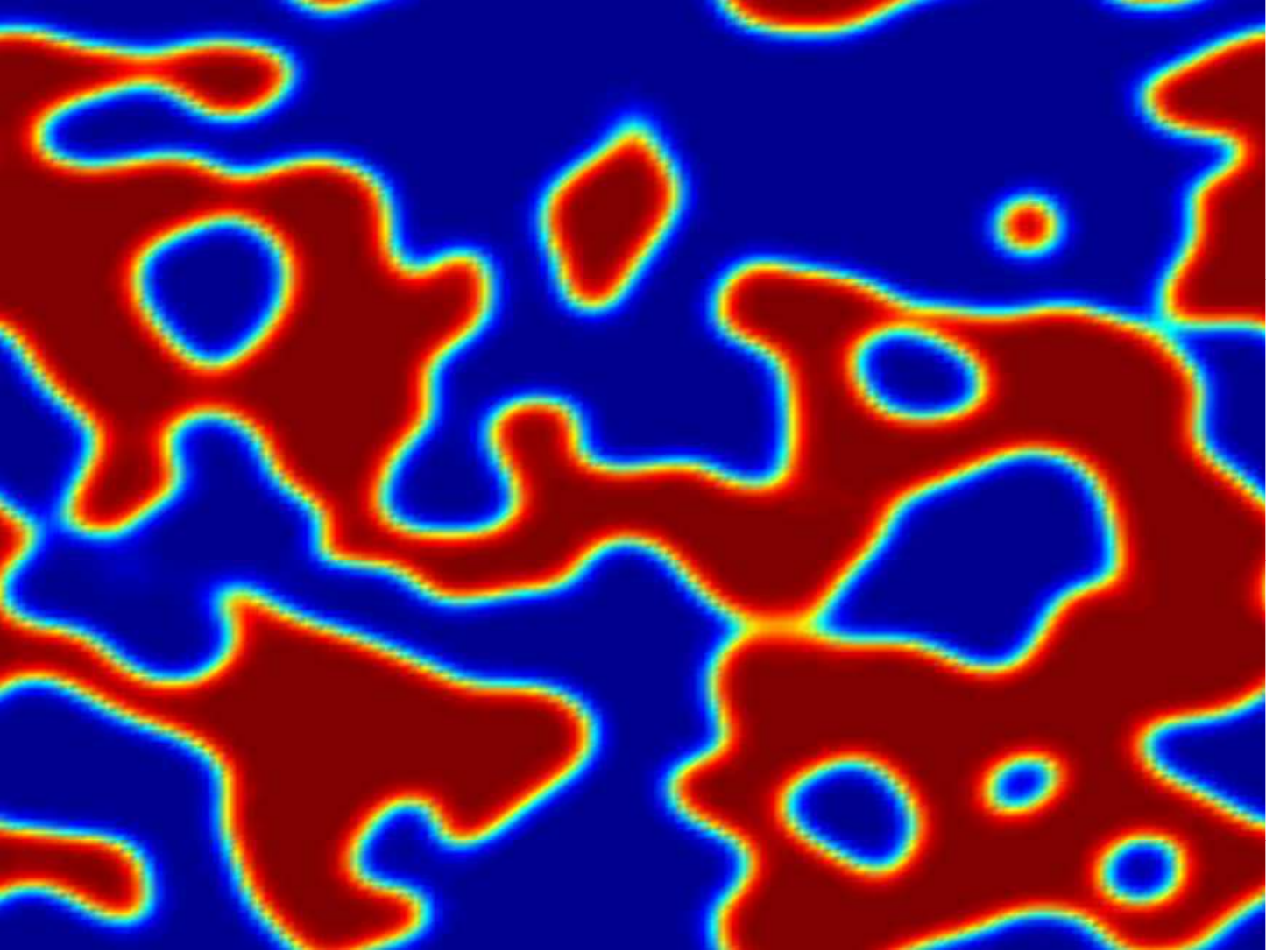}
\includegraphics[width=1.47in]{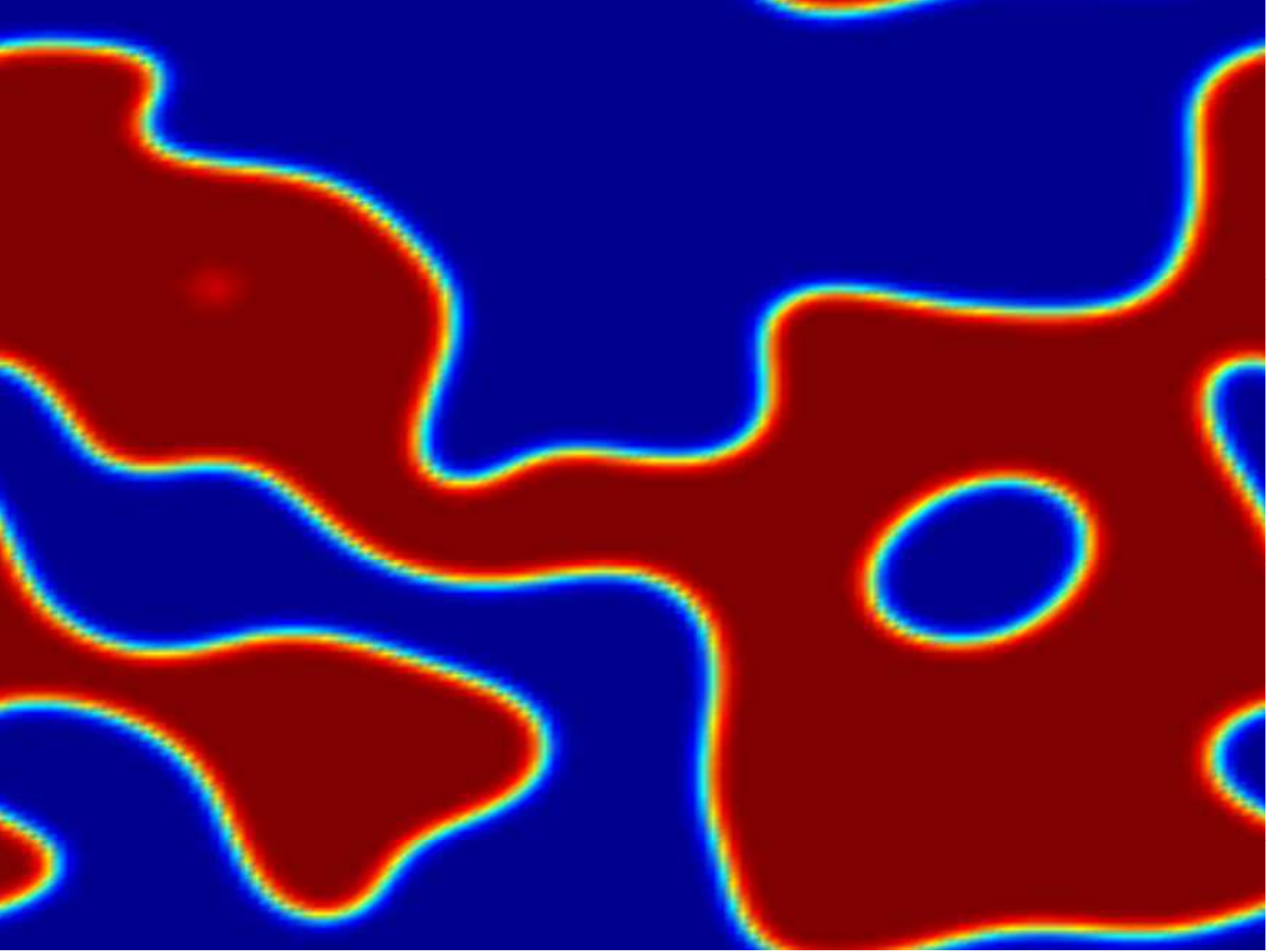}
\includegraphics[width=1.47in]{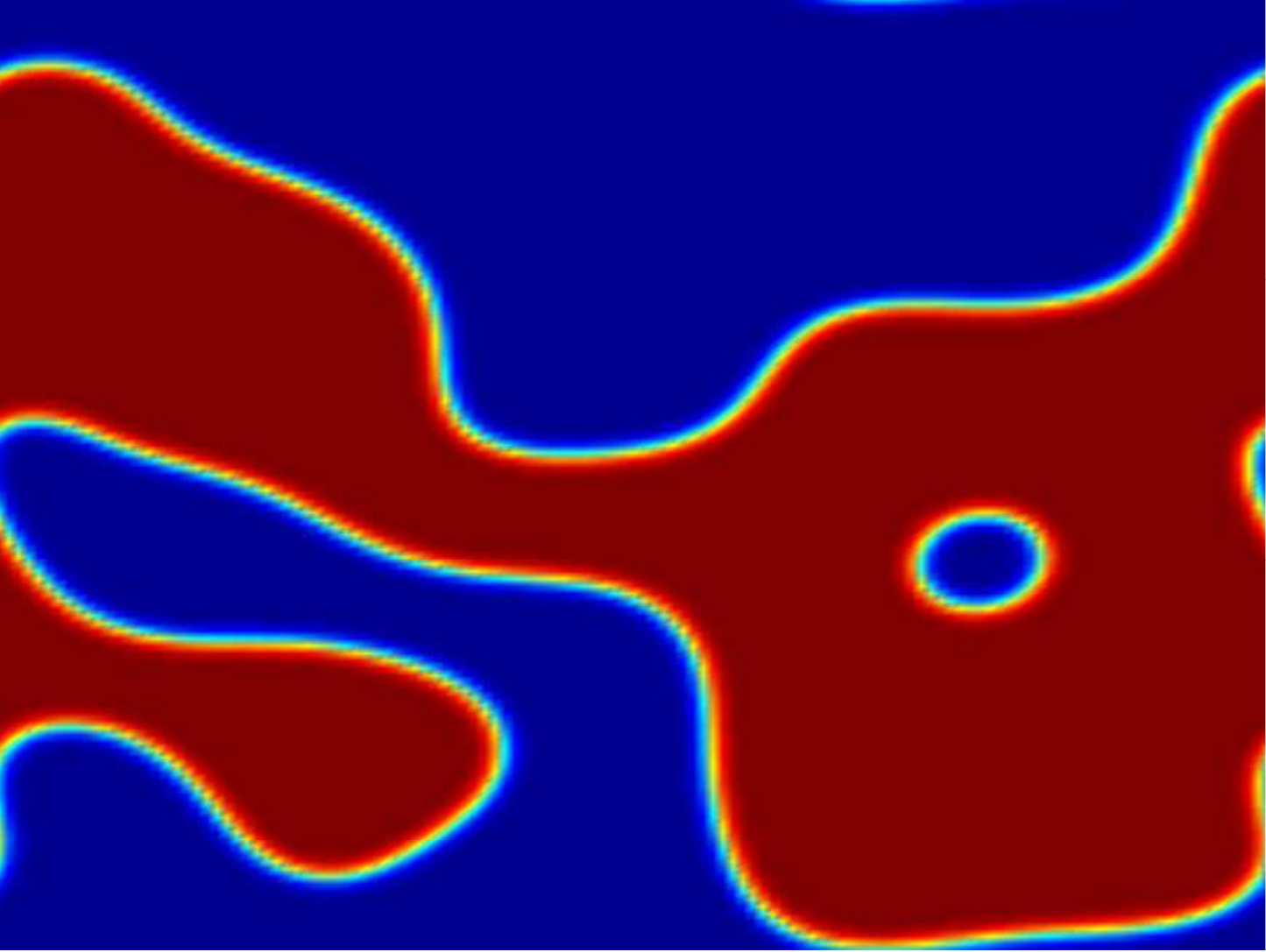}
\includegraphics[width=1.47in]{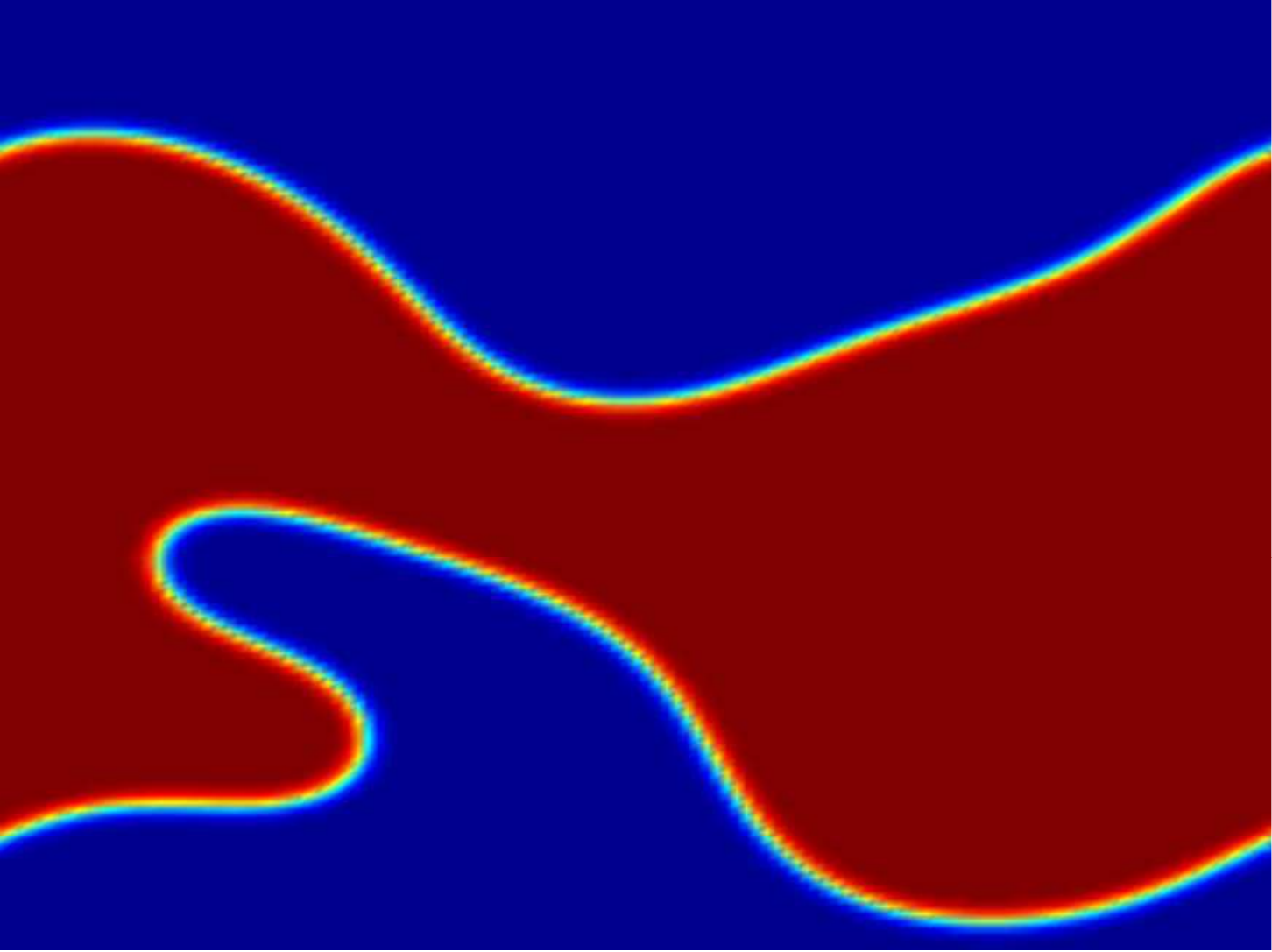}\\
\caption{Solution snapshots of Example \ref{exam:coarsen dynamic}
   at $t=10, 50, 100, 300$ (from left to right)
  for fractional orders $\alpha=0.4,\,0.7$ and $0.9$ (from top to bottom), respectively.}
\label{figs:coarsen snap}
\end{figure}

Figure \ref{figs:compar time steps} presents the discrete original energy $E(t)$,
the discrete modified energy  $E_\alpha(t)$ and different time steps
for simulating Example \ref{exam:coarsen dynamic} with $u_0(\mathbf{x})=\mathrm{rand}(\mathbf{x})$ until $T=40$.
Table \ref{table:compar CPU time} lists the CPU time and the corresponding
number of time steps for different time-stepping strategies.
The two diagrams in Figure \ref{figs:compar time steps}
show that the original and modified energies
computed on the graded-adaptive mesh
coincide with those on the graded-uniform mesh.
Table \ref{table:compar CPU time} indicates that the graded-adaptive
time-stepping strategy with appropriate parameter $\kappa$
is computationally more efficient than the graded-uniform mesh.
Also, we see that the parameter $\kappa$ affects the level of adaptivity,
i.e., the bigger the value of $\kappa$, the smaller the adaptive steps.

\begin{figure}[htb!]
\centering
\includegraphics[width=3.0in,height=2.0in]{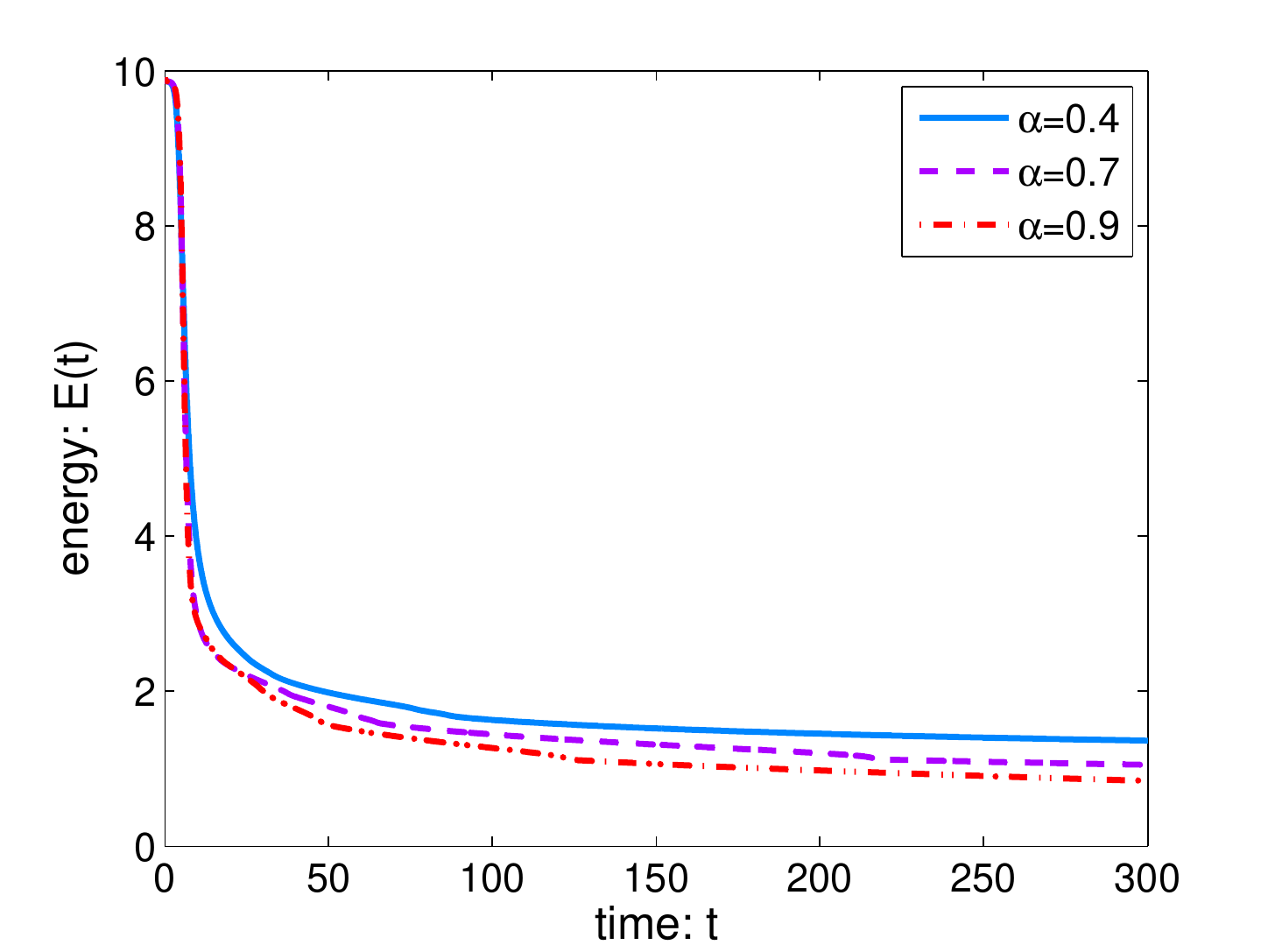}
\includegraphics[width=3.0in,height=2.0in]{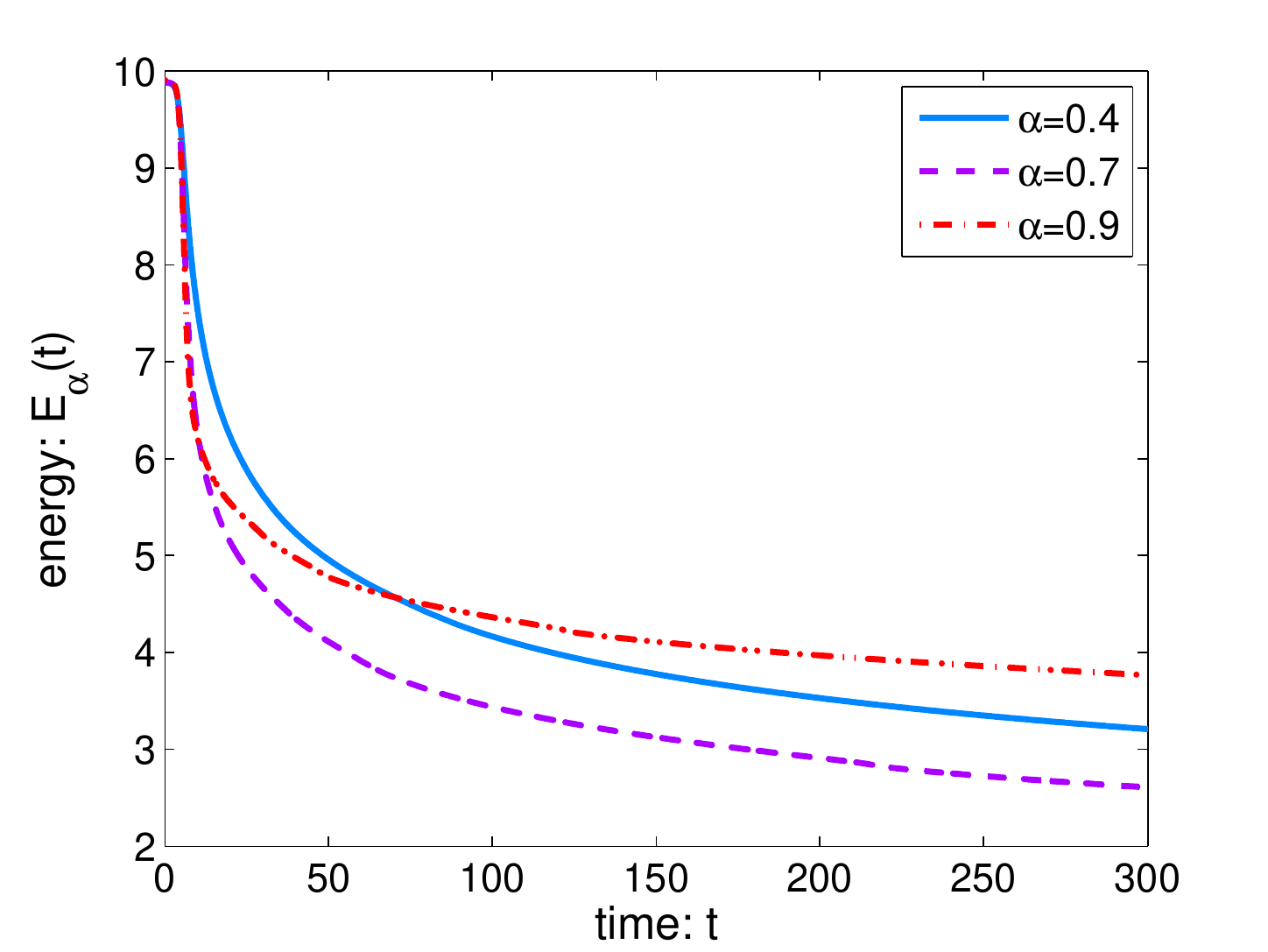}\\
\includegraphics[width=3.0in,height=2.0in]{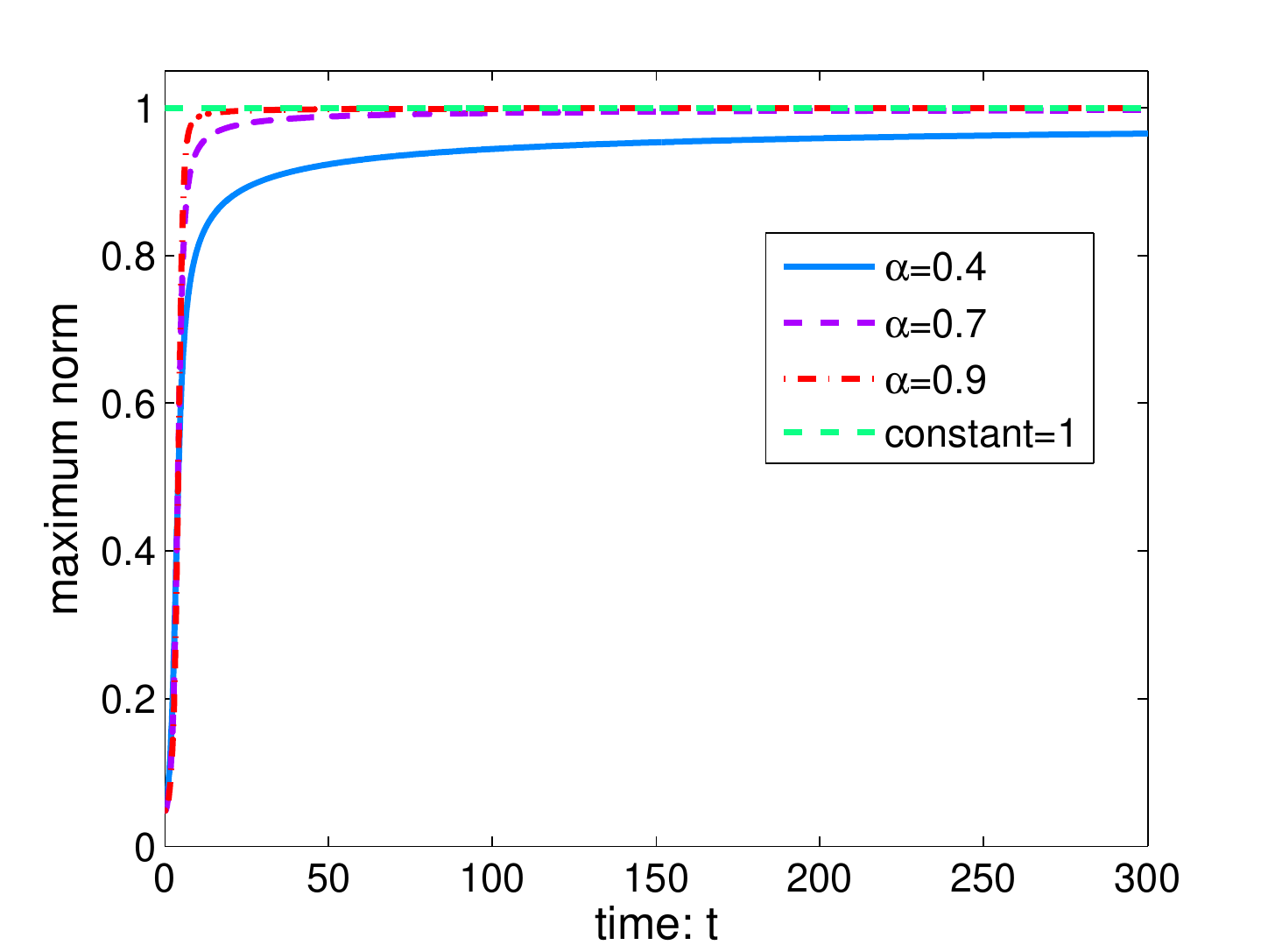}
\includegraphics[width=3.0in,height=2.0in]{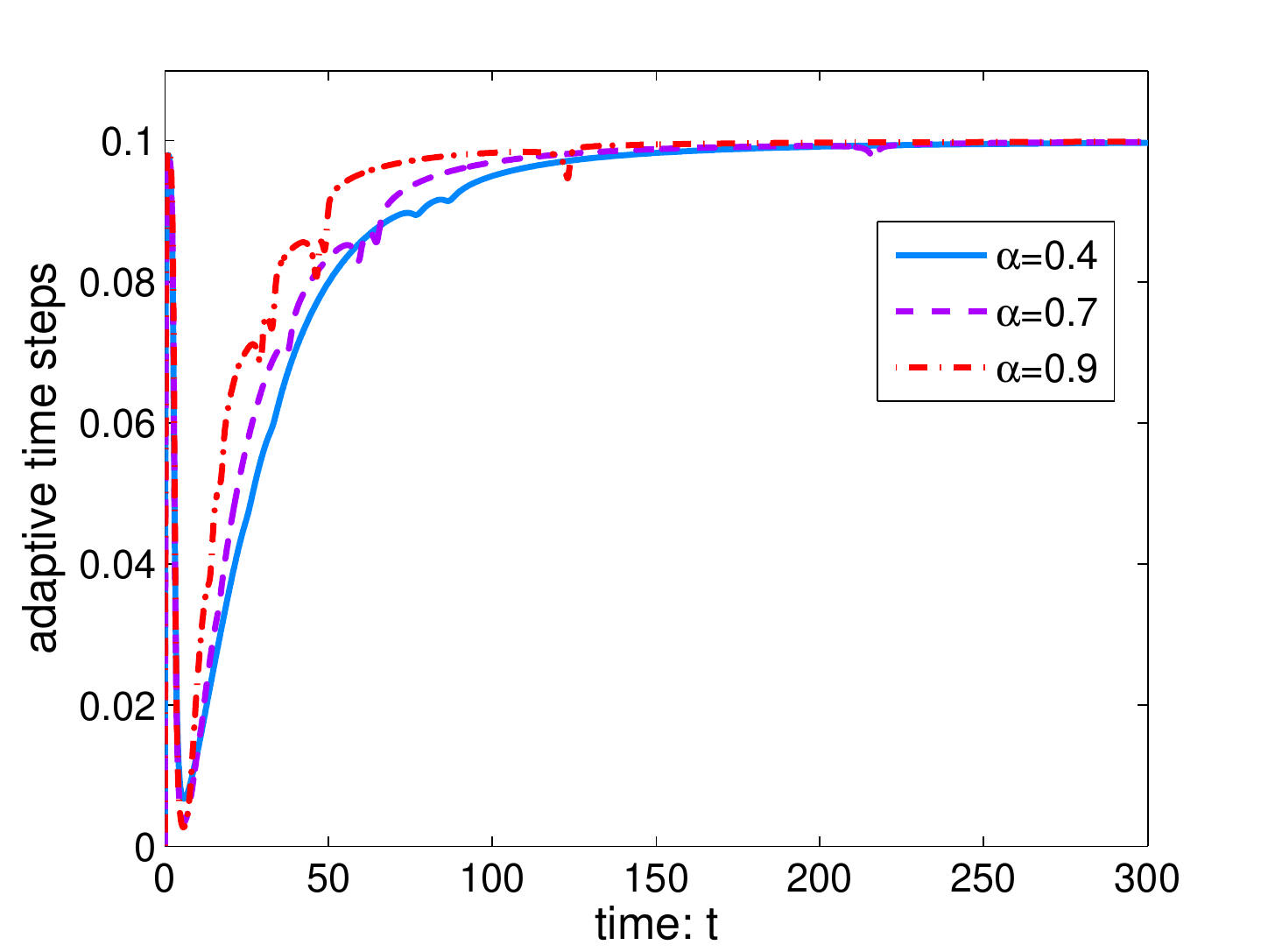}
\caption{Evolutions of energies $E(t)$ , $E_\alpha(t)$,  maximum norm and adaptive time steps
  (form left to right) of Example \ref{exam:coarsen dynamic} for three fractional orders $\alpha=0.4,\,0.7$ and $0.9$, respectively.}
\label{figs:coarsen energy max norm}
\end{figure}

Now we use the graded-adaptive time-stepping strategy with $\kappa=10^3$
to simulate the coarsening dynamics of Example \ref{exam:coarsen dynamic} until $T=300$.
The time evolutions of the microstructure due to phase separation
is summarized in Figure \ref{figs:coarsen snap}.
The time evolutions of discrete energies, $E(t)$ and $E_\alpha(t)$,
the discrete maximum principle, and adaptive time steps
are displayed in Figure \ref{figs:coarsen energy max norm}.
As seen, the coarsening dynamics for a big fractional order $\alpha$ is faster
than that for a small one. Correspondingly, the bigger the fractional order $\alpha,$
the faster the original energy $E(t)$ dissipates.
Due to the convolution term
$\mathcal{I}_{t}^{\alpha}\mynorm{\frac{\delta E}{\delta u}}^2,$
the variational energy $E_\alpha(t)$ yields a sightly different behavior.
The time-step curves show that small time  steps are selected
during the early separation progress having a large variation of energy;
large time-steps are chosen during the coarsening progress having a small variation of energy.
Moreover, the time evolution of coarsening dynamics preserves the discrete maximum principle well.

\section{Concluding remarks}
\setcounter{equation}{0}

We proposed a Crank-Nicolson scheme with variable steps for the time fractional Allen-Cahn equation that can preserve both the energy stability and the maximum bound principle. More importantly, the scheme is asymptotically energy stability preserving in the $\alpha \rightarrow 1$ limit. Our scheme is build on a reformulated problem associated with the Riemann-Liouville derivative. In this way, we build up for the first time a reversible discrete transformation between the L1-type formula of Riemann-Liouville derivative and a new L1-type formula of Caputo derivative.

This work raises some open issues to be further studied:
\begin{itemize}
  \item The numerical rate of convergence of our scheme is $1+\alpha.$ Thus it is worth to design a second order scheme that can
  preserve both the maximum bound principle and the variational energy dissipation law.
  One possible way to do this is the recently suggested second-order formula in \cite{Mustapha:2020}
  by replacing the piecewise constant approximation $\Pi_0v$
  with the piecewise linear polynomial $\Pi_1v$ in \eqref{eq: L1r formula}.

  \item By the DOC kernels \eqref{eq: orthogonal identity}, we build a connection between
  the discrete L1 Riemann-Liouville derivative \eqref{eq: L1r formula} and
  an indirect discrete Caputo derivative \eqref{eq: L1r caputo formula}.
  How about other discrete Riemann-Liouville derivatives,
  such as the variable-step second-order approximation in \cite{Mustapha:2020}?
  Lemma \ref{lem: Mutual orthogonality} suggests that there exist (indirect)
  discrete Riemann-Liouville formulas for any existing numerical Caupto derivatives,
  including the L1, Alikhanov (L2-1$_{\sigma}$) and L1$^{+}$ formulas,
  then it would be interesting to investigate  numerical approximation properties for those associated discrete Riemann-Liouville approximations.

\end{itemize}

\section*{Acknowledgements}
The authors would like to thank Dr. Bingquan Ji for his help on numerical computations.

\appendix


\section{Approximation of the nonlinear bulk}\label{append: nonlinear approximation}
\setcounter{equation}{0}

We consider a second-order approximation of the nonlinear bulk force $f(u)=u^3-u$.
It is easy to check the following equalities
\begin{align*}
4f(a)\bra{a-b}=&\,(1-a^2)^2-(1-b^2)^2-2(1-a^2)(a-b)^2
+(a^2-b^2)^2,\\
4f(b)\bra{a-b}=&\,(1-a^2)^2-(1-b^2)^2+2(1-b^2)(a-b)^2
-(a^2-b^2)^2.
\end{align*}
Then one can obtain that
\begin{align*}
\frac12\kbra{f(a)+f(b)}(a-b)=F(a)-F(b)+\frac14\brat{a+b}(a-b)^3.
\end{align*}
We consider a function $H(a,b)$ with a real parameter $\nu$,
\begin{align*}
H(a,b):=\frac12\kbra{f(a)+f(b)}-\frac{1}4\kbra{(2-\nu) a+\nu b}(a-b)^2
\end{align*}
such that
\begin{align*}
H(a,b)(a-b)=F(a)-F(b)+\frac{\nu-1}4(a-b)^4\ge F(a)-F(b)\quad\text{if $\nu\ge1$.}
\end{align*}
Moreover, the stabilized term in $H(a,b)$ contains
\begin{align*}
\kbra{(2-\nu) a+\nu b}(a-b)^2
=(2-\nu) a^3+(3\nu-4) a^2b+(2-3\nu) ab^2+\nu b^3.
\end{align*}
One can choose $\nu=4/3$ to eliminate the term $a^2b$ so that $H(a,b)$ contains only the terms $a^3$,
$ab^2$ and $b^3$. Thus we have
\begin{align}
&H(a,b)=\frac12\kbra{f(a)+f(b)}-\frac{1}{6}\brat{a+2b}(a-b)^2
=\frac{1}{3}a^{3}+\frac{1}{2}ab^{2}+\frac{1}{6}b^{3}-\frac{1}{2}(a+b),\label{fun: nonlinear accuracy}\\
&H(a,b)(a-b)=F(a)-F(b)+\frac{1}{12}(a-b)^4.\label{fun: nonlinear coercivity}
\end{align}
The function $H(a,b)$ in \eqref{fun: nonlinear accuracy} will present a second-order approximation of the function $f$ over the interval $[a,b]$.
Note that, the equality \eqref{fun: nonlinear coercivity} is vital to the unconditional
energy dissipation of the suggested method, see Theorem \ref{thm: discrete energy dissipation}.
Furthermore, the following property
\begin{align}
\frac{\partial H}{\partial a}(a,b)=a^2+\frac{1}{2}(b^2-1)\label{fun: nonlinear odd terms}
\end{align}
is important to the unique solvability and maximum bound principle of our nonlinear scheme,
see Theorem \ref{thm: unique solvability} and Theorem \ref{thm:Dis-Max-Principle}.


\section{Fast computations of L1$_{R}$ formula}\label{append: fast L1r algorithm}
\setcounter{equation}{0}

Always, the L1$_R$ fromula \eqref{eq: L1r formula} needs huge storage
and computational cost in long time simulations due to the non-locality
of Riemann-Liouville derivative \eqref{Cont: Riemann-Liouville def}.
To reduce the computational cost and memory requirements,
the sum-of-exponentials technique \cite[Theorem 2.1]{JiangZhangQianZhang:2017Fast}
is applied here to speed up the evaluation of the L1$_R$ formula.
A key result is to approximate the kernel function
$\omega_{\alpha}(t)$ efficiently inside the interval $[\Delta{t},\,T]$.

\begin{lemma}\label{lem:SOE approx}
For the given $\alpha\in(0,\,1)$, an absolute tolerance error $\epsilon\ll{1}$, a cut-off time $\Delta{t}>0$ and a finial time $T$, there exists a positive integer $N_{q}$, positive quadrature nodes $\theta^{\ell}$ and corresponding positive weights $\varpi^{\ell}\,(1\leq{\ell}\leq{N_{q}})$ such that
\begin{align*}
\bigg|
\omega_{\alpha}(t)
-\sum_{\ell=1}^{N_{q}}\varpi^{\ell}e^{-\theta^{\ell}t}\bigg|\leq\epsilon,
\quad
\forall\,{t}\in[\Delta{t},\,T],
\end{align*}
where the number $N_q$ of quadrature nodes satisfies
\[
N_q = \mathcal{O}\bra{
\log\frac{1}{\epsilon}\braB{
\log\log\frac{1}{\epsilon} + \log\frac{T}{\Delta t}}
+ \log\frac{1}{\Delta t}\braB{
\log\log\frac{1}{\epsilon} + \log\frac{1}{\Delta t}}
}.
\]
\end{lemma}

The Riemann-Liouville derivative \eqref{Cont: Riemann-Liouville def} is first split into
a local part $[t_{n-1},\,t]$ and a history part $[0,\,t_{n-1}]$.
The local part is approximated by the constant interpolation
$\bra{\Pi_{0,n}v}(t)$ and the history part is evaluated via the SOE technique, that is,
\begin{align}\label{fast:L1R approx}
\bra{{}^{R}\!\partial_{t}^{\alpha}v}(t_{n-\frac{1}{2}})
&\approx \frac{1}{\tau_n}\int_{t_{n-1}}^{t_{n}}\frac{\partial}{\partial t}
\int_{t_{n-1}}^{t}\omega_{\alpha}(t-s)(\Pi_{0,n}v)(s)\zd{s}\zd{t}
\nonumber\\
&\quad + \frac{1}{\tau_n}\int_{t_{n-1}}^{t_{n}}
\frac{\partial}{\partial t}\int_{0}^{t_{n-1}}v(s)
\sum_{\ell=1}^{N_{q}}
\varpi^{\ell}e^{-\theta^{\ell}(t-s)}\zd{s}\zd{t}\nonumber\\
&=\frac{a_{0}^{(n)}}{\tau_n}v^{n-\frac12}
+\frac{1}{\tau_n}\sum_{\ell=1}^{N_{q}}\varpi^{\ell}
\int_{t_{n-1}}^{t_{n}}\frac{\partial}{\partial t}
\int_{0}^{t_{n-1}}v(s)
e^{-\theta^{\ell}(t-s)}\zd{s}\zd{t}\nonumber\\
&=\frac{a_{0}^{(n)}}{\tau_n}v^{n-\frac12}
-\frac{1}{\tau_n}\sum_{\ell=1}^{N_{q}}\varpi^{\ell}
\brab{1-e^{-\theta^{\ell}\tau_{n}}}
\mathcal{H}^{\ell}(t_{n-1})
\quad \text{for $n\ge 1$,}
\end{align}
in which $\mathcal{H}^{\ell}(t_{k})$  is given by
\[
\mathcal{H}^{\ell}(t_{k})
:=\int_{0}^{t_{k}}e^{-\theta^{\ell}(t_{k}-s)}v(s)\zd{s}\quad
\text{with}\quad\mathcal{H}^{\ell}(t_{0})=0.
\]
Applying the constant interpolation $\Pi_{0,k}v$
to approximate $v$ in interval $[t_{k-1}, t_k]$,
one can find the following recursive formula to update $\mathcal{H}^{\ell}(t_{k})$,
\begin{align}\label{fast:history part}
\mathcal{H}^{\ell}(t_{k})
&\approx\int_{0}^{t_{k-1}}e^{-\theta^{\ell}(t_{k}-s)}v(s)\zd{s}
+\int_{t_{k-1}}^{t_{k}}e^{-\theta^{\ell}(t_{k}-s)}
v^{k-\frac12}\zd{s}\nonumber\\
&=e^{-\theta^{\ell}\tau_{k}}\mathcal{H}^{\ell}(t_{k-1})
+v^{k-\frac12}\int_{t_{k-1}}^{t_{k}}e^{-\theta^{\ell}(t_{k}-s)}\zd{s}.
\end{align}
Then the two approximations \eqref{fast:L1R approx}-\eqref{fast:history part}
gives a fast L1$_{R}$ formula,
\begin{align}\label{fast:L1R formula}
({}^{R}\!\partial_{f}^{\alpha}v)^{n-\frac{1}{2}}
:=\frac{1}{\tau_n}a_{0}^{(n)}v^{n-\frac12}
-\frac{1}{\tau_n}\sum_{\ell=1}^{N_{q}}\varpi^{\ell}
\brab{1-e^{-\theta^{\ell}\tau_{n}}}H^{\ell}(t_{n-1}),
\end{align}
where the history $\mathcal{H}^{\ell}(t_{k})$ will be updated by $H^{\ell}(t_{0})=0$ and
\begin{align*}
H^{\ell}(t_{k})
=e^{-\theta^{\ell}\tau_{k}}H^{\ell}(t_{k-1})
+v^{k-\frac12}\int_{t_{k-1}}^{t_{k}}e^{-\theta^{\ell}(t_{k}-s)}\zd{s}\quad\text{for $k\ge 1$.}
\end{align*}


\begin{thebibliography}{99}

\bibitem{AllenCahn_1979}
{\sc S. M. Allen and J. W. Cahn},
 {\em A microscopic theory for antiphase boundary motion and its
  application to antiphase domain coarsening},
 Acta Metall, 27:1085--1095, 1979.

\bibitem{AlsaediAhmadKirane:2015}
{\sc A. Alsaedi, B. Ahmad and M. Kirane},
{\em Maximum principle for certain
generalized time and space-fractional diffusion equations},
Quart. Appl. Math., 73 (2015), pp. 163-175.
%
%


\bibitem{DuYangZhou:2019}
{\sc Q.~Du, J.~Yang and Z.~Zhou},
 {\em Time-fractional {A}llen-{C}ahn equations: analysis and numerical
  methods},
 {arXiv:1906.06584v1}, 2019.


\bibitem{GomezHughes:2011Provably}
{\sc H.~Gomez and T.~J. Hughes},
 {\em Provably unconditionally stable, second-order time-accurate, mixed
  variational methods for phase-field models},
 { J. Comput. Phys.}, 230:5310--5327, 2011.

%


\bibitem{HouTangYang:2017Numerical}
{\sc T.~Hou, T.~Tang and J.~Yang},
{\em Numerical analysis of fully discretized {C}rank-{N}icolson scheme for
  fractional-in-space {Allen-Cahn} equations},
J. Sci. Comput., 72 (2017), pp. 1--18.

\bibitem{IncYusufAliyuBaleanu:2018Time}
{\sc M.~Inc, A.~Yusuf, A.~Aliyu and D.~Baleanu},
 {\em Time-fractional {A}llen-{C}ahn and time-fractional {K}lein-{G}ordon
  equations: {L}ie symmetry analysis, explicit solutions and convergence analysis},
  {Physica A Stat. Mech. Appl.}, 493:94--106, 2018.

\bibitem{JiangZhangQianZhang:2017Fast}
S. Jiang , J. Zhang, Z. Qian and Z. Zhang.
\newblock Fast evaluation of the {C}aputo fractional derivative, its
  applications to fractional diffusion equations.
\newblock {\em Comm. Comput. Phys.}, 21:650--678, 2017.


\bibitem{JiLiaoZhang:2020}
{\sc B. Ji, H.-L. Liao and L. Zhang},
{\em Simple maximum-principle preserving time-stepping methods
for time-fractional Allen-Cahn equation}, Adv. Comput. Math.,
46(2), 2020, doi: 10.1007/s10444-020-09782-2.












\bibitem{LiWangYang:2017}
{\sc Z.~Li, H.~Wang and D.~Yang},
 {\em  A space-time fractional phase-field model with tunable sharpness and
  decay behavior and its efficient numerical simulation},
 { J. Comput. Phys.}, 347:20--38, 2017.

\bibitem{LiaoLiZhang:2018}
{\sc H.-L. Liao, D. Li and J. Zhang},
{\em Sharp error estimate of nonuniform L1 formula for linear
reaction-subdiffusion equations},
SIAM J. Numer. Anal., 56(2): 1112-1133, 2018.


\bibitem{LiaoMcLeanZhang:2019}
{\sc H.-L. Liao, W. McLean and J. Zhang},
{\em A discrete {Gr\"{o}nwall} inequality
with application to numerical schemes for subdiffusion problems},
SIAM J. Numer. Anal., 57(1):218-237,  2019.


\bibitem{LiaoYanZhang:2019}
{\sc H.-L. Liao, Y. Yan and J. Zhang},
{\em Unconditional convergence of a two-level linearized fast algorithm
for semilinear subdiffusion equations},
J. Sci. Comput., 80(1):1-25, 2019.





\bibitem{LiaoTangZhou:2020jcp}
{\sc H.-L. Liao, T. Tang and T. Zhou},
{\em A second-order and nonuniform time-stepping maximum-principle preserving scheme for time-fractional
Allen-Cahn equations}, J. Comput. Phys., 414, 2020, 109473.



\bibitem{LiaoTangZhou:2020doc}
{\sc H.-L. Liao, T. Tang and T. Zhou},
{\em Positive definiteness of real quadratic forms resulting
from variable-step approximations of convolution operators},
arXiv:2011.13383v1, 2020, submitted.

\bibitem{LiaoZhang:2020}
{\sc H.-L. Liao and Z. Zhang},
{\em Analysis of adaptive BDF2 scheme for diffusion equations},
Math. Comp., 2019, DOI: 10.1090/mcom/3585.


\bibitem{LiuChengWangZhao:2018Time}
{\sc H.~Liu, A.~Cheng, H.~Wang and J.~Zhao},
 {\em  Time-fractional {A}llen-{C}ahn and {C}ahn-{H}illiard phase-field
  models and their numerical investigation},
 { Comp. Math. Appl.}, 76:1876--1892, 2018.




\bibitem{Mustapha:2011}
{\sc K. Mustapha},
\emph{An implicit finite difference time-stepping method for
a subdiffusion equation with spatial discretization by finite elements},
IMA J. Numer. Anal., 31 (2011), 719-739.

\bibitem{MustaphaAlMutawa:2012}
{\sc K. Mustapha and J. AlMutawa},
\emph{A finite difference method for an anomalous
subdiffusion equation: theory and applications},
Numer. Algor., 61 (2012), 525-543.

\bibitem{MustaphaMcLean:2009}
{\sc K. Mustapha and W. McLean},
 \emph{Piecewise-linear, discontinuous Galerkin method for a fractional
diffusion equation},
Numer. Algor., 56 (2011), 159-184.


\bibitem{Mustapha:2020}
{\sc K. Mustapha}, \emph{An L1 approximation for a fractional
reaction-diffusion equation, a second-order error
analysis over time-graded meshes},
SIAM J. Numer. Anal.,  58 (2020),  1319--1338.




\bibitem{QiaoZhangTang:2011}
{\sc Z.~Qiao, Z.~Zhang and T.~Tang},
 {\em An adaptive time-stepping strategy for the molecular beam epitaxy
  models},
 {SIAM J. Sci. Comput.}, 33:1395-1414, 2011.

\bibitem{QTY}
{\sc C. Quan, T. Tang and J. Yang},
 {\em How to define dissipation-preserving energy for time-fractional phase-field equations},
 {CSIAM Trans. Appl. Math.,} 1 (2020), pp. 478-490.












\bibitem{TangYuZhou:2018On}
{\sc T.~Tang, H.~Yu and T.~Zhou},
{\em On energy dissipation theory and numerical stability for
  time-fractional phase field equations},
SIAM J. Sci. Comput., 41-6 (2019), pp. A3757-A3778.









\bibitem{ZhaoChenWang:2019On}
{\sc J.~Zhao, L.~Chen, and H.~Wang},
{\em On power law scaling dynamics for time-fractional phase field models
  during coarsening},
  {Comm. Non. Sci. Numer. Simu.}, 70:257--270, 2019.


\end{thebibliography}
\end{document}